\documentclass[12pt,a4paper]{article}
\usepackage[titletoc,title]{appendix}
\usepackage{amsmath,amsfonts,amssymb,bm,hyperref,amsthm,tikz}
\numberwithin{equation}{section}
\newtheorem{theorem}{Theorem}[section]
\newtheorem{proposition}[theorem]{Proposition}
\newtheorem{lemma}[theorem]{Lemma}
\newtheorem{corollary}[theorem]{Corollary}
\theoremstyle{definition}

\newtheorem{remark}[theorem]{Remark}
\newtheorem{definition}[theorem]{Definition}

\usepackage[top=1in, bottom=1in, left=1in, right=1in]{geometry}
\begin{document}

\title{The second Weyl coefficient for a first order system}
\author{
Zhirayr Avetisyan
\thanks{ZA:
Department of Mathematics, University of California, Santa Barbara,
South Hall, Santa Barbara, CA 93106, USA;
Z.Avetisyan@math.ucsb.edu,
\url{http://www.z-avetisyan.com/};
this paper was written when
ZA was employed by University College London, funded by EPSRC grant EP/M000079/1.
}
\and
Johannes Sj\"ostrand\thanks{JS:
IMB, Universit\'e de Bourgogne,
9, Av.~A.~Savary, BP 47870,
FR-21780 Dijon cedex, France;
johannes.sjostrand@u-bourgogne.fr,
\url{http://sjostrand.perso.math.cnrs.fr/};
JS was supported by
EPSRC grant EP/M000079/1 and CNRS grant PRC No 1556 CNRS-RFBR 2017-2019.
}
\and
Dmitri Vassiliev\thanks{DV:
Department of Mathematics,
University College London,
Gower Street,
London WC1E~6BT,
UK;
d.vassiliev@ucl.ac.uk,
\url{http://www.ucl.ac.uk/\~ucahdva/};
DV was supported by
EPSRC grant EP/M000079/1.
}}

\renewcommand\footnotemark{}


\maketitle
\begin{abstract}
	For a scalar elliptic self-adjoint
    operator on a compact manifold without boundary we have two-term
    asymptotics for the number of eigenvalues between $0$ and
    $\lambda $ when $\lambda \to\infty $, under an additional
    dynamical condition. (See \cite[Theorem 3.5]{DuGu75} for an
    early result in this direction.) 

    \par In the case of an elliptic system of first order, the
    existence of two-term asymptotics was also established quite early
    and as in the scalar case Fourier integral operators have been the
    crucial tool. The complete computation of the coefficient of the
    second term was obtained only in the 2013 paper
    \cite{jst_part_a}. In the present paper we
    simplify that calculation. The main observation is that with the
    existence of two-term asymptotics already established, it suffices to
    study the resolvent as a pseudodifferential operator in order to
    identify and compute the second coefficient.

{\bf Keywords:} spectral theory, asymptotic distribution of eigenvalues.

\

{\bf MSC classes:} primary 35P20; secondary 35J46, 35R01.

\end{abstract}


\tableofcontents

\section{Statement of the problem}
\label{Statement of the problem}

Let $A$ be a first order linear psedodifferential operator acting on $m$-columns
of complex-valued half-densities over a connected
closed (i.e.~compact and without boundary) $n$-dimensional  mani\-fold $M$.
Throughout this paper we assume that $m,n\ge2$.

Let $A_1(x,\xi)$ and $A_\mathrm{sub}(x,\xi)$ be the principal
and subprincipal symbols of $A$.
Here $x=(x^1,\ldots,x^n)$ denotes local coordinates
and $\xi=(\xi_1,\ldots,\xi_n)$ denotes the dual variable (momentum).
The principal and subprincipal symbols are $m\times m$ matrix-functions
on $T^*M\setminus\{\xi=0\}$.

Recall that the concept of subprincipal symbol
originates from the classical paper \cite{DuiHor} of
J.~J.~Duistermaat and L.~H\"ormander: see formula (5.2.8) in that paper.
Unlike \cite{DuiHor}, we work with matrix-valued symbols, but this
does not affect the formal definition of the subprincipal symbol.

We assume our operator $A$ to be formally self-adjoint (symmetric)
with respect to the standard inner product on $m$-columns
of complex-valued half-densities,
which implies that
the principal and subprincipal symbols are Hermitian.
We also assume that our operator $A$ is elliptic:
\begin{equation}
\label{definition of ellipticity}
\det A_1(x,\xi)\ne0,\qquad\forall(x,\xi)\in T^*M\setminus\{0\}.
\end{equation}

Let $h^{(j)}(x,\xi)$ be the eigenvalues
of the matrix-function
$A_1(x,\xi)$.
Throughout this paper we assume that these are
simple for all $(x,\xi)\in T^*M\setminus\{0\}$.
The ellipticity condition
\eqref{definition of ellipticity} ensures
that all our $h^{(j)}(x,\xi)$ are nonzero.

We enumerate the eigenvalues of the principal symbol
$h^{(j)}(x,\xi)$ in increasing order, using a positive index
$j=1,\ldots,m^+$ for positive $h^{(j)}(x,\xi)$ and a negative index
$j=-1,\ldots,-m^-$ for negative $h^{(j)}(x,\xi)$. Here $m^+$ is the
number of positive eigenvalues of the principal symbol and $m^-$ is
the number of negative ones. Of course, $m^++m^-=m$.

Let $\lambda_k$ and $v_k(x)$ be
the eigenvalues and the orthonormal eigenfunctions of the operator~$A$; the
particular enumeration of these eigenvalues
(accounting for multiplicities)
is irrelevant for our purposes.
Each $v_k(x)$ is, of course,
an $m$-column of half-densities.

Let us define the two local counting functions
\begin{equation}
\label{definition of local counting functions}
N_\pm(x,\lambda):=
\begin{cases}
0\quad\text{if}\quad\lambda\le0,\\
\sum_{0<\pm\lambda_k<\lambda}\|v_k(x)\|^2\quad\text{if}\quad\lambda>0.
\end{cases}
\end{equation}
The function $N_+(x,\lambda)$ counts the eigenvalues $\lambda_k$
between zero and $\lambda$,
whereas
the function $N_-(x,\lambda)$ counts the eigenvalues $\lambda_k$
between $-\lambda$ and zero.
In both cases counting eigenvalues involves assigning them weights $\|v_k(x)\|^2$.
The quantities $\|v_k(x)\|^2$ are densities on $M$ and so are the local counting functions
$N_\pm(x,\lambda)$.

Let $\hat\rho:\mathbb{R}\to\mathbb{C}$ be a smooth function such that
$\hat\rho(t)=1$ in some neighbourhood of 0
and the support of $\hat\rho$ is sufficiently small.
Here `sufficiently small' means that
$\operatorname{supp}\hat\rho\subset(-\mathbf{T},\mathbf{T})$,
where $\mathbf{T}$ is the infimum of the lengths of all possible loops.
A loop is defined as follows.
For a given $j$, let
$(x^{(j)}(t;y,\eta),\xi^{(j)}(t;y,\eta))$ denote the Hamiltonian trajectory
originating from the point $(y,\eta)$, i.e.~solution of the system of
ordinary differential equations (the dot denotes differentiation in time $t$)
\begin{equation*}
\label{Hamiltonian system of equations}
\dot x^{(j)}=h^{(j)}_\xi(x^{(j)},\xi^{(j)}),
\qquad
\dot \xi^{(j)}=-h^{(j)}_x(x^{(j)},\xi^{(j)})
\end{equation*}
subject to the initial condition $\left.(x^{(j)},\xi^{(j)})\right|_{t=0}=(y,\eta)$.
Suppose that we have a Hamiltonian trajectory
$(x^{(j)}(t;y,\eta),\xi^{(j)}(t;y,\eta))$
and a real number $T>0$ such that
$x^{(j)}(T;y,\eta)=y$. We say in this case
that we have a loop of length $T$ originating
from the point $y\in M$.

We denote
$\rho(\lambda):=\mathcal{F}^{-1}_{t\to\lambda}[\hat\rho(t)]$,
where $\mathcal{F}^{-1}$ is the inverse Fourier transform.
See \cite[Section~6]{jst_part_a} for details.

Further on we will deal with the mollified counting functions
$(N_\pm*\rho)(x,\lambda)$ rather than the original discontinuous
counting functions $N_\pm(x,\lambda)$.
Here the star stands for convolution in the variable $\lambda$.
More specifically, we will deal with the derivative,
in the variable $\lambda$, of the mollified counting functions.
The derivative will be indicated by a prime.

It is known
\cite{
jst_review,
jst_part_a,
IvriiDoklady1980,
IvriiFuncAn1982,
ivrii_springer_lecture_notes,
ivrii_book,
kamotski,
grisha,
SafarovDSc}
that the functions $(N_\pm'*\rho)(x,\lambda)$
admit asymptotic expansions in integer powers of $\lambda\,$:
\begin{equation}
\label{expansion for mollified derivative of counting function}
(N_\pm'*\rho)(x,\lambda)=
a_{n-1}^\pm(x)\,\lambda^{n-1}
+
a_{n-2}^\pm(x)\,\lambda^{n-2}
+
a_{n-3}^\pm(x)\,\lambda^{n-3}
+\dots
\quad
\text{as}
\quad
\lambda\to+\infty.
\end{equation}

\begin{definition}
\label{definition of Weyl coefficients}
We call the coefficients
$a_k^\pm(x)$ appearing in formula
\eqref{expansion for mollified derivative of counting function}
local Weyl coefficients.
\end{definition}

Note that our definition of Weyl coefficients does not
depend on the choice of mollifier $\rho$.

It is also known
\cite{
jst_review,
jst_part_a,
IvriiDoklady1980,
IvriiFuncAn1982,
ivrii_springer_lecture_notes,
ivrii_book,
kamotski,
grisha,
SafarovDSc}
that under appropriate geometric conditions we have
\begin{equation}
\label{expansion for counting function}
N_\pm(x,\lambda)
=
\frac{a_{n-1}^\pm(x)}{n}\lambda^n
+
\frac{a_{n-2}^\pm(x)}{n-1}\lambda^{n-1}
+
o(\lambda^{n-1})
\quad
\text{as}
\quad
\lambda\to+\infty.
\end{equation}

\begin{remark}
Our Definition \ref{definition of Weyl coefficients} is somewhat nonstandard.
It is customary to call the coefficients appearing in the asymptotic
expansion \eqref{expansion for counting function} Weyl coefficients
rather than those in
\eqref{expansion for mollified derivative of counting function}.
However, for the purposes of this paper we will stick with
Definition \ref{definition of Weyl coefficients}.
\end{remark}

Further on we deal with the coefficients $a_k^+(x)$.
It is sufficient to derive formulae for the coefficients $a_k^+(x)$
because one can get formulae for $a_k^-(x)$ by replacing the operator
$A$ by the operator $-A$.

If the principal symbol of our operator $A$ is negative definite, then the operator has a finite number of
positive eigenvalues and all the coefficients $a_k^+(x)$ vanish.
So further on we assume that the principal symbol has at least one positive eigenvalue.
In other words, we assume that $m^+\ge1$.

The task at hand is to write down explicit formulae for the
coefficients $a_{n-1}^+(x)$ and $a_{n-2}^+(x)$ in terms of the principal
and subprincipal symbols of the operator $A$.

The explicit formula for the coefficient $a_{n-1}^+(x)$ has been known
since at least 1980, see, for example,
\cite{
IvriiDoklady1980,
IvriiFuncAn1982,
ivrii_springer_lecture_notes,
ivrii_book,
kamotski,
grisha,
SafarovDSc}.
It reads
\begin{equation}
\label{formula for a plus}
a_{n-1}^+(x)=\frac{n}{(2\pi)^n}\,\sum_{j=1}^{m^+}
\ \int\limits_{h^{(j)}(x,\xi)<1}d\xi\,,
\end{equation}
where $d\xi=d\xi_1\ldots d\xi_n$.

The explicit formula for the coefficient $a_{n-2}^+(x)$ was derived only in 2013, see
\cite[formula (1.24)]{jst_part_a}. This formula reads
\begin{multline}
\label{formula for b plus}
a_{n-2}^+(x)=-\frac{n(n-1)}{(2\pi)^n}\,\sum_{j=1}^{m^+}
\ \int\limits_{h^{(j)}(x,\xi)<1}
\Bigl(
[v^{(j)}]^*A_\mathrm{sub}v^{(j)}
\\
-\frac i2
\{
[v^{(j)}]^*,A_1-h^{(j)},v^{(j)}
\}
+\frac i{n-1}h^{(j)}\{[v^{(j)}]^*,v^{(j)}\}
\Bigr)(x,\xi)\,
d\xi\,.
\end{multline}
Here curly brackets denote the Poisson bracket on matrix-functions
$\{P,R\}:=P_{x^\alpha}R_{\xi_\alpha}-P_{\xi_\alpha}R_{x^\alpha}$
and its further generalisation
\begin{equation}
\label{definition of generalised Poisson bracket}
\{F,G,H\}:=F_{x^\alpha}GH_{\xi_\alpha}-F_{\xi_\alpha}GH_{x^\alpha}\,,
\end{equation}
where the subscripts $x^\alpha$ and $\xi_\alpha$
indicate partial derivatives and
the repeated index $\alpha$ indicates summation over $\alpha=1,\ldots,n$.

Note that if $q(x,\xi)$ is a function on $T^*M\setminus\{0\}$ positively
homogeneous in $\xi$ of degree 0, then
\[
\int\limits_{h^{(j)}(x,\xi)<1}q(x,\xi)\,d\xi
\]
is a density on $M$.
Hence, the quantities
\eqref{formula for a plus}
and
\eqref{formula for b plus}
are densities.

The problem with the derivation of formula \eqref{formula for b plus}
given in \cite{jst_part_a}
was that it was very complicated.
The aim of the current paper is to provide an alternative,
much simpler, derivation of formula \eqref{formula for b plus}.

It may be that the approach outlined in the current paper would allow one,
in the future, to calculate further coefficients in the asymptotic
expansion
\eqref{expansion for mollified derivative of counting function}.
Note that for an operator that is not semibounded this is a nontrivial task.

\section{Strategy for the evaluation of the second Weyl coefficient}
\label{Strategy for the evaluation of the second Weyl coefficient}

Let $z\in\mathbb{C}$, $\operatorname{Im}z>0$.
Our basic idea is to consider the resolvent
$(A-zI)^{-1}$ and, by studying it, recover the second Weyl coefficient $a_{n-2}^+(x)$.
Unfortunately, the operator $(A-zI)^{-1}$ is not of trace class, therefore one has to
modify our basic idea
so as to reduce our analysis to that of trace class operators.

Let us consider the self-adjoint operator
\begin{equation}
\label{Strategy formula 1}
i\left[
2(A-zI)^{1-n}-(A-2zI)^{1-n}
-2(A-\bar zI)^{1-n}+(A-2\bar zI)^{1-n}
\right].
\end{equation}
We claim that the operator \eqref{Strategy formula 1} is of trace class.
In order to justify this claim we calculate below, for fixed $z$, the principal symbol
of the operator \eqref{Strategy formula 1} and
show that it has degree of homogeneity $-n-1$.

Let $B$ be the parametrix (approximate inverse) of $A$,
see \cite[Section 5]{shubin} for details.
Then, modulo $L^{-\infty}(M)$
(integral operators with infinitely smooth integral kernels),
we have
\[
A-zI\equiv A-zAB=A(I-zB),
\]
\[
(A-zI)^{n-1}\equiv A^{n-1}(I-zB)^{n-1},
\]
\begin{equation}
\label{expansion 1}
(A-zI)^{1-n}\equiv (I-zB)^{1-n}A^{1-n}\equiv (I-zB)^{1-n}B^{n-1}.
\end{equation}
But
\begin{equation}
\label{expansion 2}
(I-zB)^{1-n}\equiv I+(n-1)zB-\frac{n(n-1)}2(zB)^2+\ldots,
\end{equation}
where the expansion is understood as an asymptotic expansion in smoothness
(each subsequent term is a pseudodifferential operator of lower order).
Substituting
\eqref{expansion 2}
into
\eqref{expansion 1},
we get
\begin{equation}
\label{expansion 3}
(A-zI)^{1-n}
\equiv B^{n-1}
+(n-1)zB^n-\frac{n(n-1)}2z^2B^{n+1}
+\ldots.
\end{equation}
Replacing $z$ by $2z$, we get
\begin{equation}
\label{expansion 4}
(A-2zI)^{1-n}
\equiv B^{n-1}
+2(n-1)zB^n-2n(n-1)z^2B^{n+1}
+\ldots.
\end{equation}
Formulae
\eqref{expansion 3}
and
\eqref{expansion 4}
imply
\begin{equation}
\label{expansion 5}
2(A-zI)^{1-n}
-
(A-2zI)^{1-n}
\equiv B^{n-1}
+
n(n-1)z^2B^{n+1}
+\ldots.
\end{equation}
Replacing $z$ by $\bar z$, we get
\begin{equation}
\label{expansion 6}
2(A-\bar zI)^{1-n}
-
(A-2\bar zI)^{1-n}
\equiv B^{n-1}
+
n(n-1)\bar z^2B^{n+1}
+\ldots.
\end{equation}
Formulae
\eqref{expansion 5}
and
\eqref{expansion 6}
imply that the operator
\eqref{Strategy formula 1}
is a pseudodifferential operator of order $-n-1$
with principal symbol
$-4n(n-1)(\operatorname{Re}z)(\operatorname{Im}z)A_1^{-n-1}$.

It might seem more natural to consider the operator
\begin{equation}
\label{Strategy formula 1a}
(A-zI)^{-n-1}
\end{equation}
instead of \eqref{Strategy formula 1}.
The operator \eqref{Strategy formula 1a}
is also of order $-n-1$, hence, trace class.
Unfortunately, the algorithm presented in the remainder
of this section won't work for
the operator \eqref{Strategy formula 1a}.
The reason is that if
we start with
\eqref{Strategy formula 1a},
we end up with the integral
\begin{equation}
\label{Strategy formula 1b}
\int_0^{+\infty}\frac{\mu^{n-2}}{(\mu-z)^{n+1}}\,d\mu\,,
\end{equation}
where the exponent in the numerator is lower that the exponent in the denominator
by more than one.
The integral
\eqref{Strategy formula 1b}
is a polynomial in $\,\frac1z\,$
(no logarithm!) and it does not experience a jump when
$z$ crosses the positive real axis.
Starting with
\eqref{Strategy formula 1a}
one can recover
$a_{n-2}^+-(-1)^na_{n-2}^-\,$,
but it appears to be impossible to recover
$a_{n-2}^+$ itself.
We need a logarithm in order to separate contributions
from positive and negative eigenvalues.

The operator
\eqref{Strategy formula 1}
is a pseudodifferential operator of order $\,-n-1\,$,
hence it has a continuous integral kernel.
This observation allows us to introduce the following definition.

\begin{definition}
\label{definition of density f}
By $f(x,z)$ we denote the real-valued continuous density
obtained by restricting the integral kernel
of the operator \eqref{Strategy formula 1}
to the diagonal $x=y$ and taking the matrix trace~$\,\operatorname{tr}\,$.
\end{definition}

The explicit formula for our density is
\begin{equation}
\label{Strategy formula 2}
f(x,z)=
i
\sum_{\lambda_k}
\left[
\frac{2}{(\lambda_k-z)^{n-1}}
-
\frac{1}{(\lambda_k-2z)^{n-1}}
-
\frac{2}{(\lambda_k-\bar z)^{n-1}}
+
\frac{1}{(\lambda_k-2\bar z)^{n-1}}
\right]
\|v_k(x)\|^2\,.
\end{equation}
This formula can be equivalently rewritten as
\begin{multline}
\label{Strategy formula 3}
f(x,z)
=
i
\int_0^{+\infty}
\left[
\frac{2}{(\mu-z)^{n-1}}
-
\frac{1}{(\mu-2z)^{n-1}}
-
\frac{2}{(\mu-\bar z)^{n-1}}
+
\frac{1}{(\mu-2\bar z)^{n-1}}
\right]
N_+'(x,\mu)\,d\mu
\\
-
(-1)^n
\frac{2^n-1}{2^{n-1}}
i
\left[
\frac1{z^{n-1}}
-
\frac1{\bar z^{n-1}}
\right]
\sum_{\lambda_k=0}
\|v_k(x)\|^2
\\
-
(-1)^n
\,i
\int_0^{+\infty}
\left[
\frac{2}{(\mu+z)^{n-1}}
-
\frac{1}{(\mu+2z)^{n-1}}
-
\frac{2}{(\mu+\bar z)^{n-1}}
+
\frac{1}{(\mu+2\bar z)^{n-1}}
\right]
N_-'(x,\mu)\,d\mu\,.
\end{multline}
The expression in the second line of \eqref{Strategy formula 3} is the contribution
from the kernel (eigenspace corresponding to the eigenvalue zero)
of the operator $A$.

Let us also introduce another density
\begin{multline}
\label{Strategy formula 4}
f^\rho(x,z):=
i
\int_0^{+\infty}
\left[
\frac{2}{(\mu-z)^{n-1}}
-
\frac{1}{(\mu-2z)^{n-1}}
-
\frac{2}{(\mu-\bar z)^{n-1}}
+
\frac{1}{(\mu-2\bar z)^{n-1}}
\right]
(N_+'*\rho)(x,\mu)\,d\mu
\\
-
(-1)^n
\,i
\int_0^{+\infty}
\left[
\frac{2}{(\mu+z)^{n-1}}
-
\frac{1}{(\mu+2z)^{n-1}}
-
\frac{2}{(\mu+\bar z)^{n-1}}
+
\frac{1}{(\mu+2\bar z)^{n-1}}
\right]
(N_-'*\rho)(x,\mu)\,d\mu\,.
\end{multline}

Put $z=\lambda e^{i\varphi}$, where $\lambda>0$
and $0<\varphi<\pi$.
We will now fix the angle $\varphi$ and examine what happens when $\lambda\to+\infty$.

\begin{lemma}
\label{Strategy lemma 1}
The density $f^\rho(x,\lambda e^{i\varphi})-f(x,\lambda e^{i\varphi})$ tends to zero
as $\lambda\to+\infty$.
\end{lemma}

\textbf{Proof\ }
See Appendix \ref{Proof of Lemma Strategy lemma 1}.
\qed

\begin{lemma}
\label{Strategy lemma 2}
The density $f^\rho(x,\lambda e^{i\varphi})$
admits the asymptotic expansion
\begin{equation}
\label{Strategy formula 5}
f^\rho(x,\lambda e^{i\varphi})
=
b_1(x,\varphi)
\lambda
+
b_0(x,\varphi)
+o(1)
\quad\text{as}\quad
\lambda\to+\infty,
\end{equation}
where
\begin{equation}
\label{Strategy formula 6}
b_1(x,\varphi)=
-4(\ln 2)(n-1)(\sin\varphi)
\left[
a_{n-1}^+(x)+(-1)^n\,a_{n-1}^-(x)
\right],
\end{equation}
\begin{equation}
\label{Strategy formula 7}
b_0(x,\varphi)=
-2
\left[
(\pi-\varphi)\,a_{n-2}^+(x)+(-1)^n\,\varphi\,a_{n-2}^-(x)
\right].
\end{equation}
\end{lemma}

\textbf{Proof\ }
See Appendices \ref{Some integrals involving the functions} and \ref{Proof of Lemma Strategy lemma 2}.
\qed

Lemmata
\ref{Strategy lemma 1}
and
\ref{Strategy lemma 2}
imply the following corollary.

\begin{corollary}
The density $f(x,\lambda e^{i\varphi})$
admits the asymptotic expansion
\begin{equation}
\label{Strategy formula 8}
f(x,\lambda e^{i\varphi})
=
b_1(x,\varphi)
\lambda
+
b_0(x,\varphi)
+o(1)
\quad\text{as}\quad
\lambda\to+\infty,
\end{equation}
where the coefficients $b_1(x,\varphi)$ and $b_0(x,\varphi)$
are given by formulae
\eqref{Strategy formula 6}
and
\eqref{Strategy formula 7}
respectively.
\end{corollary}

Suppose that we know the coefficient  $b_0(x,\varphi)$ for all $\varphi\in(0,\pi)$.
It is easy to see that formula
\eqref{Strategy formula 7}
allows us to recover the second Weyl coefficient $a_{n-2}^+(x)$.
Namely, if we take an arbitrary pair of distinct $\varphi_1,\varphi_2\in(0,\pi)$ then
\begin{equation}
\label{Strategy formula 9}
a_{n-2}^+(x)=
\frac{\varphi_1\,b_0(x,\varphi_2)-\varphi_2\,b_0(x,\varphi_1)}{2\pi(\varphi_2-\varphi_1)}\,.
\end{equation}
Alternatively, the second Weyl coefficient $a_{n-2}^+(x)$ can be recovered by means of the identity
\begin{equation}
\label{Strategy formula 10}
a_{n-2}^+(x)=
-\frac1{2\pi}
\lim_{\varphi\to0^+}
b_0(x,\varphi)\,.
\end{equation}

Formulae
\eqref{Strategy formula 8}--\eqref{Strategy formula 10}
tell us that
the problem of evaluating the second Weyl coefficient
has been reduced to evaluating the second coefficient in the asymptotic
expansion of the density $f(x,\lambda e^{i\varphi})$ as $\lambda\to+\infty$.
Recall that the latter is defined in accordance with
Definition \ref{definition of density f}.

\section{The Weyl symbol of the resolvent}
\label{The Weyl symbol of the resolvent}

Let $z=\lambda e^{i\varphi}$, where $\lambda>0$
and $0<\varphi<\pi$.
We formally assign to $z$ a `weight',
as if it were positively homogeneous in $\xi$ of degree 1.
Our argument goes along the lines of
\cite[Section 9]{shubin}.

We performed formal calculations evaluating the symbol of the operator
$(A-zI)^{-1}$  in local coordinates
and then switched to the Weyl symbol.
(One could have worked with Weyl symbols from the very start.)
Further on we denote the Weyl symbol of the operator
$(A-zI)^{-1}$ by $[(A-zI)^{-1}]_W$.
We calculated $[(A-zI)^{-1}]_W$ in the two leading terms:
\begin{multline}
\label{12 December 2016 equation 1}
[(A-zI)^{-1}]_W
=
(A_1-zI)^{-1}
-
(A_1-zI)^{-1}A_\mathrm{sub}(A_1-zI)^{-1}
\\
+\frac i2
\{
(A_1-zI)^{-1},A_1-zI,(A_1-zI)^{-1}
\}
+O[(1+|\xi|+|z|)^{-2}(1+|\xi|)^{-1}].
\end{multline}
Here the curly brackets denote the generalised Poisson
bracket on matrix functions \eqref{definition of generalised Poisson bracket}.

The concept of a Weyl symbol was initially introduced for
pseudodifferential operators in $\mathbb{R}^n$, see \cite[subsection 23.3]{shubin}.
In the case of pseudodifferential operators acting on half-densities
over a manifold it turns out that the Weyl symbol depends on the choice of
local coordinates.
However, in the two leading terms the Weyl symbol does not
depend on the choice of local coordinates,
see Appendix \ref{Weyl quantization on manifolds}.
Note that a consistent definition of the full Weyl symbol
for a pseudodifferential operator acting on half-densities
over a manifold requires the introduction of an affine connection, see
\cite{mckeag}. In the current paper we do not assume that we have a connection.

See
Appendix \ref{Symbolic approximation for the resolvent and its powers}
for a discussion of symbol classes
and an explanation of the origins of the particular structure of
the remainder term in formula \eqref{12 December 2016 equation 1},
as well as remainder term estimates in subsequent
formulae.
In (\ref{E.21.5}) we obtain (\ref{12 December 2016 equation 1}) in the
appropriate symbol classes.

Note that the expression in the second line of
\eqref{12 December 2016 equation 1}
can be equivalently rewritten as
\begin{equation}
\label{12 December 2016 equation 2}
\{
(A_1-zI)^{-1},A_1-zI,(A_1-zI)^{-1}
\}
=
(A_1-zI)^{-1}
\{
A_1,(A_1-zI)^{-1},A_1
\}
(A_1-zI)^{-1},
\end{equation}
which is the representation used by V.~Ivrii,
see second displayed formula on page 226 of~\cite{ivrii_springer_lecture_notes}.
We mention \eqref{12 December 2016 equation 2} in order to put our analysis within
the context of previous research in the subject.

Let us now express the
principal symbol $A_1$ in terms of its
eigenvalues $h^{(j)}$ and eigenprojections $P^{(j)}$:
\begin{equation}
\label{12 December 2016 equation 3}
A_1=\sum_jh^{(j)}P^{(j)}.
\end{equation}
In what follows we will be substituting
\eqref{12 December 2016 equation 3} into
our previous formulae.
But before proceeding with the calculations let us discuss
which expression,
the one in the RHS of \eqref{12 December 2016 equation 2}
or
the one in the LHS of \eqref{12 December 2016 equation 2},
is better suited for practical purposes.
Substitution of \eqref{12 December 2016 equation 3} into the
RHS of \eqref{12 December 2016 equation 2} gives a sum over five indices,
whereas
substitution of \eqref{12 December 2016 equation 3} into the
LHS of \eqref{12 December 2016 equation 2} gives a sum over only three indices.
Hence, we will stick with the representation from the
LHS of \eqref{12 December 2016 equation 2}.

Substituting
\eqref{12 December 2016 equation 3}
into \eqref{12 December 2016 equation 1}
we get
\begin{multline}
\label{Weyl symbol of resolvent decomposed 1}
[(A-zI)^{-1}]_W=
\sum_j\frac{P^{(j)}}{h^{(j)}-z}
-
\sum_{k,l}\frac{P^{(k)}A_\mathrm{sub}P^{(l)}}{(h^{(k)}-z)(h^{(l)}-z)}
\\
+\frac i2
\sum_{j,k,l}
(h^{(j)}-z)
\left\{
\frac{P^{(k)}}{h^{(k)}-z}
\,,
P^{(j)}
,
\frac{P^{(l)}}{h^{(l)}-z}
\right\}
+O[(1+|\xi|+|z|)^{-2}(1+|\xi|)^{-1}].
\end{multline}
Our eigenprojections satisfy the identity
\begin{equation}
\label{basic identity for eigenprojections}
P^{(k)}P^{(j)}=\delta^{kj}P^{(k)}.
\end{equation}
The identity
\eqref{basic identity for eigenprojections}
allows us to rewrite formula
\eqref{Weyl symbol of resolvent decomposed 1}
as
\begin{multline}
\label{Weyl symbol of resolvent decomposed 2}
[(A-zI)^{-1}]_W=
\sum_j\frac{P^{(j)}}{h^{(j)}-z}
-
\sum_{k,l}\frac{P^{(k)}A_\mathrm{sub}P^{(l)}}{(h^{(k)}-z)(h^{(l)}-z)}
\\
+\frac i2
\sum_{j,k,l}
\frac{h^{(j)}-z}
{(h^{(k)}-z)(h^{(l)}-z)}
\{
P^{(k)},P^{(j)},P^{(l)}
\}
\\
-\frac i2
\sum_{k,l}
\frac
{
P^{(k)}
\bigl(
h^{(k)}_{x^\alpha}
P^{(l)}_{\xi_\alpha}
-
h^{(k)}_{\xi_\alpha}
P^{(l)}_{x^\alpha}
\bigr)
+
\bigl(
h^{(l)}_{\xi_\alpha}
P^{(k)}_{x^\alpha}
-
h^{(l)}_{x^\alpha}
P^{(k)}_{\xi_\alpha}
\bigr)
P^{(l)}
}
{(h^{(k)}-z)(h^{(l)}-z)}
\\
+O[(1+|\xi|+|z|)^{-2}(1+|\xi|)^{-1}].
\end{multline}

\section{The matrix trace of the resolvent}
\label{The matrix trace of the resolvent}

Let $B$ be a matrix pseudodifferential operator
acting on $m$-columns of half-densities, $v\mapsto Bv$.
The action of such an operator can be written in more detailed form as
\begin{equation}
\label{The matrix trace of the resolvent equation 1}
\begin{pmatrix}
v_1\\
v_2\\
\vdots\\
v_m
\end{pmatrix}
\mapsto
\begin{pmatrix}
B_1{}^1&B_1{}^2&\dots&B_1{}^m\\
B_2{}^1&B_2{}^2&\dots&B_2{}^m\\
\vdots&\vdots&\ddots&\vdots\\
B_m{}^1&B_m{}^2&\dots&B_m{}^m
\end{pmatrix}
\begin{pmatrix}
v_1\\
v_2\\
\vdots\\
v_m
\end{pmatrix},
\end{equation}
where the $B_j{}^k$ are scalar pseudodifferential operators
acting on half-densities.

\begin{definition}
The matrix trace of the operator
\eqref{The matrix trace of the resolvent equation 1}
is the scalar operator
\begin{equation}
\label{The matrix trace of the resolvent equation 2}
\operatorname{tr}B:=
B_1{}^1+B_2{}^2+\dots+B_m{}^m.
\end{equation}
\end{definition}

Obviously, the Weyl symbol of the matrix trace of an operator
is the matrix trace of the Weyl symbol of the operator.
Hence, formula
\eqref{Weyl symbol of resolvent decomposed 2} implies
\begin{multline}
\label{The matrix trace of the resolvent equation 3}
[\operatorname{tr}(A-zI)^{-1}]_W=
\sum_j\frac{1}{h^{(j)}-z}
-
\sum_j\frac{\operatorname{tr}[A_\mathrm{sub}P^{(j)}]}{(h^{(j)}-z)^2}
\\
+\frac i2
\sum_{j,k,l}
\frac{h^{(j)}-z}
{(h^{(k)}-z)(h^{(l)}-z)}
\operatorname{tr}
\{
P^{(k)},P^{(j)},P^{(l)}
\}
+O[(1+|\xi|+|z|)^{-2}(1+|\xi|)^{-1}].
\end{multline}
Note that formula \eqref{The matrix trace of the resolvent equation 3}
does not contain terms with derivatives of the Hamiltonians $h^{(j)}$
because all such terms cancelled out after we took the matrix trace.

Formula \eqref{basic identity for eigenprojections} implies
\begin{multline}
\label{20 December 2016 equation 1}
\operatorname{tr}\{P^{(k)},P^{(j)},P^{(l)}\}
=2\delta^{kj}\delta^{jl}\operatorname{tr}\{P^{(j)},P^{(j)},P^{(j)}\}
\\
-\delta^{kj}\operatorname{tr}\{P^{(l)},P^{(j)},P^{(l)}\}
-\delta^{jl}\operatorname{tr}\{P^{(k)},P^{(j)},P^{(k)}\}
+\delta^{kl}\operatorname{tr}\{P^{(k)},P^{(j)},P^{(k)}\}.
\end{multline}
Substituting
\eqref{20 December 2016 equation 1}
into
\eqref{The matrix trace of the resolvent equation 3}
and using
\eqref{12 December 2016 equation 3}
we get
\begin{multline}
\label{20 December 2016 equation 2}
[\operatorname{tr}(A-zI)^{-1}]_W=
\sum_j\frac1{h^{(j)}-z}
-
\sum_j\frac{\operatorname{tr}[A_\mathrm{sub}P^{(j)}]}{(h^{(j)}-z)^2}
\\
+\frac i2
\sum_j
\frac{
\operatorname{tr}
\{
P^{(j)},A_1-h^{(j)}I,P^{(j)}
\}
}
{(h^{(j)}-z)^2}
+i
\sum_j
\frac{
\operatorname{tr}
\{
P^{(j)},P^{(j)},P^{(j)}
\}
}
{h^{(j)}-z}
\\
+O[(1+|\xi|+|z|)^{-2}(1+|\xi|)^{-1}].
\end{multline}
Detailed calculations leading up to formulae
\eqref{20 December 2016 equation 1}
and
\eqref{20 December 2016 equation 2}
are presented in Appendix~\ref{Proof of formulae (4.4) and (4.5)}.

Formula \eqref{20 December 2016 equation 2}
provides a compact representation for the
Weyl symbol of the matrix trace of the resolvent.
Even though our intermediate calculations involved summation
over several (up to three) indices,
summation in our final formula \eqref{20 December 2016 equation 2}
is carried out over a single index.

\section{The matrix trace of a power of the resolvent}
\label{The matrix trace of a power of the resolvent}

In order to implement the strategy outlined in
Section \ref{Strategy for the evaluation of the second Weyl coefficient}
we need to write down the Weyl symbol of the operator
$\,\operatorname{tr}(A-zI)^{1-n}\,$.

We have the operator identity
\begin{equation}
\label{The matrix trace of a power of the resolvent equation 1}
(A-zI)^{1-n}=
\frac1{(n-2)!}
\,
\frac{d^{n-2}}{dz^{n-2}}
\,
(A-zI)^{-1}
\,.
\end{equation}
The operations of taking the matrix trace
and differentiation with respect to a parameter commute,
so formula \eqref{The matrix trace of a power of the resolvent equation 1}
implies
\begin{equation}
\label{The matrix trace of a power of the resolvent equation 2}
\operatorname{tr}
(A-zI)^{1-n}=
\frac1{(n-2)!}
\,
\frac{d^{n-2}}{dz^{n-2}}
\,
\operatorname{tr}(A-zI)^{-1}
\,.
\end{equation}
The latter formula, in turn, implies
\begin{equation}
\label{The matrix trace of a power of the resolvent equation 3}
[\operatorname{tr}
(A-zI)^{1-n}]_W=
\frac1{(n-2)!}
\,
\frac{d^{n-2}}{dz^{n-2}}
\,
[\operatorname{tr}(A-zI)^{-1}]_W
\,.
\end{equation}

Substituting
\eqref{20 December 2016 equation 2}
into
\eqref{The matrix trace of a power of the resolvent equation 3}
we get
\begin{multline}
\label{The matrix trace of a power of the resolvent equation 4}
[\operatorname{tr}(A-zI)^{1-n}]_W=
\sum_j\frac1{(h^{(j)}-z)^{n-1}}
-
(n-1)
\sum_j\frac{\operatorname{tr}[A_\mathrm{sub}P^{(j)}]}{(h^{(j)}-z)^n}
\\
+\frac i2
(n-1)
\sum_j
\frac{
\operatorname{tr}
\{
P^{(j)},A_1-h^{(j)}I,P^{(j)}
\}
}
{(h^{(j)}-z)^n}
+i
\sum_j
\frac{
\operatorname{tr}
\{
P^{(j)},P^{(j)},P^{(j)}
\}
}
{(h^{(j)}-z)^{n-1}}
\\
+O[(1+|\xi|+|z|)^{-n}(1+|\xi|)^{-1}].
\end{multline}
We can view this as an explicit version of the result of
applying $(d/dz)^{n-2}$ to the trace of (\ref{E.21.5}) (cf.\ (\ref{E.37})).

\section{Asymptotic expansion for the density $f$}
\label{Asymptotic expansion for the density}

We have previously defined the density $f(x,z)$, see
Definition \ref{definition of density f}.
In this section we shall derive the asymptotic expansion for
the density $f(x,\lambda e^{i\varphi})$ as $\lambda\to+\infty$.
The angle $0<\varphi<\pi$ will be assumed to be fixed.

Put
\begin{equation}
\label{Asymptotic expansion for the density equation 1}
s_{1-n}^{(j)}(x,\xi,z):=
\frac1{(h^{(j)}-z)^{n-1}}\,,
\end{equation}
\begin{multline}
\label{Asymptotic expansion for the density equation 2}
s_{-n}^{(j)}(x,\xi,z):=
-
(n-1)
\frac{\operatorname{tr}[A_\mathrm{sub}P^{(j)}]}{(h^{(j)}-z)^n}
+\frac i2
(n-1)
\frac{
\operatorname{tr}
\{
P^{(j)},A_1-h^{(j)}I,P^{(j)}
\}
}
{(h^{(j)}-z)^n}
\\
+i
\frac{
\operatorname{tr}
\{
P^{(j)},P^{(j)},P^{(j)}
\}
}
{(h^{(j)}-z)^{n-1}}\,,
\end{multline}
where the subscripts indicate the degree of homogeneity in $\xi$.
Recall, yet again, that our convention is `$z$ and $\xi$ are of the same order'.
Comparing
\eqref{The matrix trace of a power of the resolvent equation 4}
with
\eqref{Asymptotic expansion for the density equation 1}
and
\eqref{Asymptotic expansion for the density equation 2}
we see that
$\sum_js_{1-n}^{(j)}$ is the leading (principal) component of the Weyl symbol of the operator
$\operatorname{tr}(A-zI)^{1-n}$,
whereas
$\sum_js_{-n}^{(j)}$ is the next (subprincipal) component.

The structure of formula
\eqref{Asymptotic expansion for the density equation 1} is very simple,
whereas the structure of formula
\eqref{Asymptotic expansion for the density equation 2} is nontrivial.
This warrants a discussion.

The first term in the RHS of \eqref{Asymptotic expansion for the density equation 2}
contains the expression $\operatorname{tr}[A_\mathrm{sub}P^{(j)}]$.
It gives the `obvious' contribution to the second Weyl coefficient.
The expression $\operatorname{tr}[A_\mathrm{sub}P^{(j)}]$
appears in the early papers of V.~Ivrii and G.~V.~Rozenblyum.

The second term in the RHS of \eqref{Asymptotic expansion for the density equation 2}
contains the expression 
$\operatorname{tr}
\{
P^{(j)},A_1-h^{(j)}I,P^{(j)}
\}
$.
It gives a contribution to the second Weyl coefficient which is not so obvious.
The expression
$\operatorname{tr}
\{
P^{(j)},A_1-h^{(j)}I,P^{(j)}
\}
$
first appeared in \cite{SafarovDSc}.

Finally, the third term in the RHS of \eqref{Asymptotic expansion for the density equation 2}
contains the expression 
$\operatorname{tr}
\{
P^{(j)},P^{(j)},P^{(j)}
\}
$.
It gives a $\mathrm{U}(1)$ curvature contribution to the second Weyl coefficient.
This contribution to the second Weyl coefficient
was identified in \cite{jst_part_a} and did not appear in
previous publications.

The density $f(x,\lambda e^{i\varphi})$ is the value of the integral
kernel of the operator
\begin{equation}
\label{Asymptotic expansion for the density equation 3}
i\operatorname{tr}
\left[
2(A-zI)^{1-n}-(A-2zI)^{1-n}
-2(A-\bar zI)^{1-n}+(A-2\bar zI)^{1-n}
\right]
\end{equation}
on the diagonal.
We obtain the asymptotic expansion
\eqref{Strategy formula 8}
for $f(x,\lambda e^{i\varphi})$
by replacing the operator \eqref{Asymptotic expansion for the density equation 3}
with its Weyl symbol and integrating in $\xi$.
This gives the following formulae for the asymptotic coefficients:
\begin{equation}
\label{Asymptotic expansion for the density equation 4}
b_1(x,\varphi)
=\frac1{(2\pi)^n}
\sum_j b_1^{(j)}(x,\varphi),
\end{equation}
\begin{equation}
\label{Asymptotic expansion for the density equation 5}
b_0(x,\varphi)
=\frac1{(2\pi)^n}
\sum_j b_0^{(j)}(x,\varphi),
\end{equation}
where
\begin{multline}
\label{Asymptotic expansion for the density equation 6}
b_1^{(j)}(x,\varphi)
=
\\
i\int\left[
2s_{1-n}^{(j)}(x,\xi,e^{i\varphi})
-
s_{1-n}^{(j)}(x,\xi,2e^{i\varphi})
-
2s_{1-n}^{(j)}(x,\xi,e^{-i\varphi})
+
s_{1-n}^{(j)}(x,\xi,2e^{-i\varphi})
\right]
d\xi,
\end{multline}
\begin{multline}
\label{Asymptotic expansion for the density equation 7}
b_0^{(j)}(x,\varphi)
=
\\
i\int\left[
2s_{-n}^{(j)}(x,\xi,e^{i\varphi})
-
s_{-n}^{(j)}(x,\xi,2e^{i\varphi})
-
2s_{-n}^{(j)}(x,\xi,e^{-i\varphi})
+
s_{-n}^{(j)}(x,\xi,2e^{-i\varphi})
\right]
d\xi.
\end{multline}
The integrands in
\eqref{Asymptotic expansion for the density equation 6}
and
\eqref{Asymptotic expansion for the density equation 7}
decay as $|\xi|^{-n-1}$ as $|\xi|\to+\infty$,
so these integrals converge.

Strictly speaking, we also have to consider the
contributions from the terms $K^{(n)}$ in (\ref{E.33}). However, it follows
from the remark after (\ref{E.35}) that they are $o(1)$ as $\lambda \to+\infty$.

\section{The second Weyl coefficient}
\label{The second Weyl coefficient}

Let us us examine what happens to the integral
\eqref{Asymptotic expansion for the density equation 7}
when $\varphi\to0^+$.
It is easy to see that if $j$ is such that $h^{(j)}<0$ then the integral
\eqref{Asymptotic expansion for the density equation 7}
tends to zero as $\varphi\to0^+$: one can simply set $\varphi=0$
in the integrand.
This means that only those $j$ for which $h^{(j)}>0$ contribute to the limit of
the expression \eqref{Asymptotic expansion for the density equation 6}
when $\varphi\to0^+$.
Therefore, formulae
\eqref{Strategy formula 10}
and
\eqref{Asymptotic expansion for the density equation 5}
give us the following expression for the second Weyl coefficient:
\begin{equation}
\label{The second Weyl coefficient equation 1}
a_{n-2}^+(x)=
-\frac1{(2\pi)^{n+1}}
\sum_{j=1}^{m^+}
\lim_{\varphi\to0^+}
b_0^{(j)}(x,\varphi)\,.
\end{equation}
Here the enumeration of eigenvalues of the principal symbol $A_1$
is assumed to be chosen in such a way that
$j=1,\dots,m^+$ correspond to positive eigenvalues $h^{(j)}$.

It remains only to evaluate
$\,\lim_{\varphi\to0^+}
b_0^{(j)}(x,\varphi)\,$
explicitly.
Here $b_0^{(j)}(x,\varphi)$ is defined by formula
\eqref{Asymptotic expansion for the density equation 7},
where the integrand is defined in accordance with
\eqref{Asymptotic expansion for the density equation 2}.

Let us rewrite formula
\eqref{Asymptotic expansion for the density equation 2}
as
\begin{equation}
\label{The second Weyl coefficient equation 2}
s_{-n}^{(j)}(x,\xi,z)=s_{-n}^{(j;1)}(x,\xi,z)+s_{-n}^{(j;2)}(x,\xi,z),
\end{equation}
where
\begin{equation}
\label{The second Weyl coefficient equation 3}
s_{-n}^{(j;1)}(x,\xi,z):=
-
(n-1)
\frac{
\operatorname{tr}
\left(
A_\mathrm{sub}P^{(j)}
-\frac i2
\{
P^{(j)},A_1-h^{(j)}I,P^{(j)}
\}
\right)}
{(h^{(j)}-z)^n}
\,,
\end{equation}
\begin{equation}
\label{The second Weyl coefficient equation 4}
s_{-n}^{(j;2)}(x,\xi,z):=
i
\frac{
h^{(j)}
\operatorname{tr}
\{
P^{(j)},P^{(j)},P^{(j)}
\}
}
{h^{(j)}(h^{(j)}-z)^{n-1}}\,.
\end{equation}
Note that the numerators
in
\eqref{The second Weyl coefficient equation 3}
and
\eqref{The second Weyl coefficient equation 4}
are positively homogeneous in $\xi$ of degree zero.

Formula
\eqref{Asymptotic expansion for the density equation 7}
now reads
\begin{equation}
\label{The second Weyl coefficient equation 5}
b_0^{(j)}(x,\varphi)
=
b_0^{(j;1)}(x,\varphi)
+
b_0^{(j;2)}(x,\varphi),
\end{equation}
where
\begin{multline}
\label{The second Weyl coefficient equation 6}
b_0^{(j;k)}(x,\varphi)
=
\\
i\int\left[
2s_{-n}^{(j;k)}(x,\xi,e^{i\varphi})
-
s_{-n}^{(j;k)}(x,\xi,2e^{i\varphi})
-
2s_{-n}^{(j;k)}(x,\xi,e^{-i\varphi})
+
s_{-n}^{(j;k)}(x,\xi,2e^{-i\varphi})
\right]
d\xi\,,
\end{multline}
$k=1,2$.

Denote by $(S^*_xM)^{(j)}$ the $(n-1)$-dimensional unit cosphere in the cotangent fibre
defined by the equation $h^{(j)}(x,\xi)=1$
and denote by $d(S^*_xM)^{(j)}$ the
surface area element on $(S^*_xM)^{(j)}$
defined by the condition
\begin{equation}
\label{The second Weyl coefficient equation 7}
\left[
\frac{d}{d\mu}
\int_{h^{(j)}(x,\xi)<\mu}g(\xi)\,d\xi
\right]_{\mu=1}
=
\int_{(S^*_xM)^{(j)}}
g(\xi)
\,
d(S^*_xM)^{(j)}\,,
\end{equation}
where $g:\mathbb{R}^n\to\mathbb{R}$ is an arbitrary smooth function.
This means that we introduce spherical coordinates in the cotangent fibre
with the Hamiltonian $h^{(j)}$
playing the role of the radial coordinate, see also \cite[subsection 1.1.10]{mybook}.

Switching to spherical coordinates, we see that
each integral \eqref{The second Weyl coefficient equation 6}
is a product of two integrals,
an $(n-1)$-dimensional surface integral over the unit cosphere
and a 1-dimensional integral over the radial coordinate.
Namely, we have
\begin{equation}
\label{The second Weyl coefficient equation 8}
b_0^{(j;k)}(x,\varphi)
=
c^{(j;k)}(x)\,
d^{(j;k)}(\varphi)\,,
\end{equation}
where
\begin{equation}
\label{The second Weyl coefficient equation 9}
c^{(j;1)}(x):=
-(n-1)
\int_{(S^*_xM)^{(j)}}
\operatorname{tr}
\left(
A_\mathrm{sub}P^{(j)}
-\frac i2
\{
P^{(j)},A_1-h^{(j)}I,P^{(j)}
\}
\right)
d(S^*_xM)^{(j)}\,,
\end{equation}
\begin{equation}
\label{The second Weyl coefficient equation 10}
c^{(j;2)}(x):=
i
\int_{(S^*_xM)^{(j)}}
h^{(j)}
\operatorname{tr}
\{
P^{(j)},P^{(j)},P^{(j)}
\}
\,
d(S^*_xM)^{(j)}\,,
\end{equation}
\begin{equation}
\label{The second Weyl coefficient equation 11}
d^{(j;1)}(\varphi):=
i
\int_0^{+\infty}
\left[
\frac{2}{(\mu-e^{i\varphi})^n}
-
\frac{1}{(\mu-2e^{i\varphi})^n}
-
\frac{2}{(\mu-e^{-i\varphi})^n}
+
\frac{1}{(\mu-2e^{-i\varphi})^n}
\right]
\mu^{n-1}\,d\mu
\,,
\end{equation}
\begin{multline}
\label{The second Weyl coefficient equation 12}
d^{(j;2)}(\varphi):=
\\
i
\int_0^{+\infty}
\left[
\frac{2}{(\mu-e^{i\varphi})^{n-1}}
-
\frac{1}{(\mu-2e^{i\varphi})^{n-1}}
-
\frac{2}{(\mu-e^{-i\varphi})^{n-1}}
+
\frac{1}{(\mu-2e^{-i\varphi})^{n-1}}
\right]
\mu^{n-2}\,d\mu
\,.
\end{multline}

Integrating by parts we see that the integrals in the right--hand-sides of
\eqref{The second Weyl coefficient equation 11}
and
\eqref{The second Weyl coefficient equation 12}
have the same values,
i.e.~they do not depend on $n$.
Hence, it is sufficient to evaluate the integral
\eqref{The second Weyl coefficient equation 12}
for $n=2$.
We have
\begin{multline}
\label{The second Weyl coefficient equation 13}
d^{(j;1)}(\varphi)=
d^{(j;2)}(\varphi)
=
i
\int_0^{+\infty}
\left[
\frac{2}{\mu-e^{i\varphi}}
-
\frac{1}{\mu-2e^{i\varphi}}
-
\frac{2}{\mu-e^{-i\varphi}}
+
\frac{1}{\mu-2e^{-i\varphi}}
\right]
d\mu
\\
=-2(\pi-\varphi)
\,,
\end{multline}
so substituting
\eqref{The second Weyl coefficient equation 5},
\eqref{The second Weyl coefficient equation 8}
and
\eqref{The second Weyl coefficient equation 13}
into
\eqref{The second Weyl coefficient equation 1}
we get
\begin{equation}
\label{The second Weyl coefficient equation 14}
a_{n-2}^+(x)=
\frac1{(2\pi)^n}
\sum_{j=1}^{m^+}
\bigl[
c^{(j;1)}(x)
+
c^{(j;2)}(x)
\bigr]\,.
\end{equation}

Formulae
\eqref{The second Weyl coefficient equation 14},
\eqref{The second Weyl coefficient equation 9}
and
\eqref{The second Weyl coefficient equation 10}
give us the required explicit representation of the second Weyl coefficient.
However, integrating over a unit cosphere is not very convenient, so we rewrite
formulae
\eqref{The second Weyl coefficient equation 9}
and
\eqref{The second Weyl coefficient equation 10}
as
\begin{equation}
\label{The second Weyl coefficient equation 15}
c^{(j;1)}(x)=
-n(n-1)
\int_{h^{(j)}(x,\xi)<1}
\operatorname{tr}
\left(
A_\mathrm{sub}P^{(j)}
-\frac i2
\{
P^{(j)},A_1-h^{(j)}I,P^{(j)}
\}
\right)
(x,\xi)\,
d\xi\,,
\end{equation}
\begin{equation}
\label{The second Weyl coefficient equation 16}
c^{(j;2)}(x)=
n\,i
\int_{h^{(j)}(x,\xi)<1}
\left(
h^{(j)}\,
\operatorname{tr}
\{
P^{(j)},P^{(j)},P^{(j)}
\}
\right)
(x,\xi)\,
d\xi\,.
\end{equation}

Working with eigenprojections $P^{(j)}$ is also not very convenient, so we express them
via the normalised eigenvectors $v^{(j)}$
of the principal symbol $A_1$ as
\begin{equation}
\label{The second Weyl coefficient equation 17}
P^{(j)}=v^{(j)}[v^{(j)}]^*.
\end{equation}
Substituting
\eqref{The second Weyl coefficient equation 17}
into
\eqref{The second Weyl coefficient equation 15}
and
\eqref{The second Weyl coefficient equation 16}
we get
\begin{equation}
\label{The second Weyl coefficient equation 18}
c^{(j;1)}(x)=
-n(n-1)
\int_{h^{(j)}(x,\xi)<1}
\left(
[v^{(j)}]^*A_\mathrm{sub}v^{(j)}
-\frac i2
\{
[v^{(j)}]^*,A_1-h^{(j)}I,v^{(j)}
\}
\right)
(x,\xi)\,
d\xi\,,
\end{equation}
\begin{equation}
\label{The second Weyl coefficient equation 19}
c^{(j;2)}(x)=
-n\,i
\int_{h^{(j)}(x,\xi)<1}
\left(
h^{(j)}\,
\{
[v^{(j)}]^*,v^{(j)}
\}
\right)
(x,\xi)\,
d\xi\,.
\end{equation}
The transition from
\eqref{The second Weyl coefficient equation 15}
to
\eqref{The second Weyl coefficient equation 18}
is quite straightforward,
but the transition from
\eqref{The second Weyl coefficient equation 16}
to
\eqref{The second Weyl coefficient equation 19}
warrants an explanation.
Here we have
$\operatorname{tr}
\{
P^{(j)},P^{(j)},P^{(j)}
\}
=
-
\operatorname{tr}
(
P^{(j)}
\{
P^{(j)},P^{(j)}
\}
)
=
-
\{
[v^{(j)}]^*,v^{(j)}
\}
$,
where at the last step we made use of \cite[formula (4.17)]{jst_part_a}.

The advantage of formulae
\eqref{The second Weyl coefficient equation 18}
and
\eqref{The second Weyl coefficient equation 19}
is that they do not involve the matrix trace.

Combining formulae
\eqref{The second Weyl coefficient equation 14},
\eqref{The second Weyl coefficient equation 18}
and
\eqref{The second Weyl coefficient equation 19}
we arrive at
\eqref{formula for b plus}.

\section{Acknowledgements}
The authors are grateful to Yan-Long Fang for valuable suggestions.

\begin{appendices}

\section{Proof of Lemma \ref{Strategy lemma 1}}
\label{Proof of Lemma Strategy lemma 1}

Let us introduce the functions
\begin{equation}
\label{definition of g}
g_n(\mu,z):=\frac2{(\mu-z)^n}-\frac1{(\mu-2z)^n}-\mbox{c.c.},\quad n\in\mathbb{N},
\quad\mu\in\mathbb{R},
\quad z\in\mathbb{C}\setminus\mathbb{R}.
\end{equation}
Here and further on `$\mbox{c.c.}$' stands for `complex conjugate terms'.

The functions \eqref{definition of g} possess the following properties:
\begin{equation}
\label{BPartDerivProp}
\partial_1 g_n(\mu,z):=\partial_\mu g_n(\mu,z)=-ng_{n+1}(\mu,z),
\end{equation}
\begin{eqnarray}
\label{BBound}
|g_n(\mu,z)|\le\frac4{|\mu-z|^n}+\frac2{|\mu-2z|^n}\,.
\end{eqnarray}

Formula \eqref{Strategy formula 4}
can be rewritten as
\begin{eqnarray}\label{frhodef}
f^\rho(x,z)=i\int_0^{+\infty}g_{n-1}(\mu,z)\,(N_+'*\rho)(x,\mu)\,d\mu\nonumber\\
-(-1)^ni\int_0^{+\infty}g_{n-1}(\mu,-z)\,(N_-'*\rho)(x,\mu)\,d\mu,
\end{eqnarray}
where
\begin{equation}\label{NpmDef}
N_\pm'(x,\nu)=\sum_{\pm\lambda_k>0}\delta(\nu\mp\lambda_k)\|v_k(x)\|^2
\end{equation}
is a tempered distribution in $\nu$ supported on $\mathbb{R}_+$ and taking values in densities. The convolution
\begin{equation}
\label{dima1}
(N_\pm'*\rho)(x,\mu)=\int_0^{+\infty}N_\pm'(x,\nu)\,\rho(\mu-\nu)\,d\nu
\end{equation}
is a continuous function of $\mu$ taking values in densities. It is known that
\begin{equation*}
|(N_\pm'*\rho)(x,\mu)|\le c(x)(1+|\mu|^{n-1}),
\end{equation*}
where $c(x)$ is a fixed positive density.
Arguing as in
\eqref{expansion 1}--\eqref{expansion 6},
it is easy to see that, for fixed~$z$,
the function $g_{n-1}(\mu,z)$ decays as $|\mu|^{-n-1}$
when $\mu\to\pm\infty$,
so the integrals in (\ref{frhodef}) converge.

We have
\begin{multline}
\label{dima2}
\int_0^{+\infty}g_{n-1}(\mu,z)\,(N_\pm'*\rho)(x,\mu)\,d\mu
\\
=
\int_0^{+\infty}g_{n-1}(\mu,z)
\left(
\int_0^{+\infty}N_\pm'(x,\nu)\,\rho(\mu-\nu)\,d\nu
\right)
d\mu
\\
=
\int_0^{+\infty}
N_\pm'(x,\nu)
\left(
\int_0^{+\infty}
g_{n-1}(\mu,z)\,\rho(\mu-\nu)\,d\mu
\right)
d\nu
\\
=\int_0^{+\infty}N_\pm'(x,\mu)
\left(
\int_0^{+\infty}
g_{n-1}(\nu,z)\,\rho(\nu-\mu)\,d\nu
\right)
d\mu.
\end{multline}
In going from the second line of \eqref{dima2}
to the third
we changed the order of integration.
This can be justified, for example, by replacing
the infinite series \eqref{NpmDef}
by a finite partial sum
and going to the limit.

Substituting \eqref{dima2}
into (\ref{frhodef}) and using formula
\eqref{Strategy formula 3},
we find that
\begin{multline*}
f^\rho(x,z)-f(x,z)
=i\int_0^{+\infty}
N_+'(x,\mu)\left(
\int_0^{+\infty}
g_{n-1}(\nu,z)\,\rho(\nu-\mu)\,d\nu-g_{n-1}(\mu,z)\right)d\mu
\\
-(-1)^ni\int_0^{+\infty}
N_-'(x,\mu)\left(
\int_0^{+\infty}
g_{n-1}(\nu,-z)\,\rho(\nu-\mu)\,d\nu-g_{n-1}(\mu,-z)\right)d\mu
\\
+(-1)^n\frac{2^n-1}{2^{n-1}}i\left[\frac1{z^{n-1}}-\frac1{\bar z^{n-1}}\right]
\sum_{\lambda_k=0}\|v_k(x)\|^2.
\end{multline*}

Now, let $z=\lambda e^{i\varphi}$ with $\lambda>0$
and fixed $\varphi\in(0,\pi)$.
In view of the fact that
$N_\pm(x,\lambda)=O(\lambda^n)$, in order to show that
$f^\rho(x,\lambda e^{i\varphi})-f(x,\lambda e^{i\varphi})\to0$
as $\lambda\to+\infty$ it is sufficient to prove that
\begin{equation}\label{ToProveThat}
\left|\int_0^{+\infty}
g_{n-1}(\nu,\lambda e^{i\varphi})\,\rho(\nu-\mu)\,d\nu
-g_{n-1}(\mu,\lambda e^{i\varphi})\right|
\le\frac{\operatorname{const}_\varphi}{\lambda(1+\mu^{n+1})}\,,
\quad\forall\lambda\ge1,\quad\forall\mu\ge0.
\end{equation}

Recall that according to our definition of the mollifier $\rho$ we have
\begin{equation}\label{rhoSchwartz}
|\rho(\nu)|\le\frac{c_p}{(1+|\nu|)^p}\,,\quad\forall p\in\mathbb{N},
\end{equation}
\begin{equation}\label{rhoProp}
\int_{-\infty}^{+\infty}\rho(\nu)d\nu=1,\quad\mbox{and}\quad\int_{-\infty}^{+\infty}\rho(\nu)\nu^m d\nu=0,\quad\forall m\in\mathbb{N}.
\end{equation}

Formula (\ref{rhoProp}) implies that
\begin{multline}\label{Estimand}
\int_0^{+\infty}g_{n-1}(\nu,\lambda e^{i\varphi})\rho(\nu-\mu)d\nu-g_{n-1}(\mu,\lambda e^{i\varphi})\\
=\int_{-\infty}^{+\infty}\left[g_{n-1}(\nu,\lambda e^{i\varphi})
-g_{n-1}(\mu,\lambda e^{i\varphi})\right]\rho(\nu-\mu)d\nu-\int_{-\infty}^0g_{n-1}(\nu,\lambda e^{i\varphi})\rho(\nu-\mu)d\nu.
\end{multline}
Using (\ref{BBound}) and (\ref{rhoSchwartz}) with $p=n+3$ we get
\begin{eqnarray}\label{Est1}
\left|\int_{-\infty}^0g_{n-1}(\nu,\lambda e^{i\varphi})\rho(\nu-\mu)d\nu\right|\le\int_{-\infty}^0\frac6{\lambda^{n-1}|\sin\varphi|^{n-1}}\frac{c_{n+3}}{(1+|\nu|+\mu)^{n+3}}d\nu\nonumber\\
\le\frac{6c_{n+3}}{\lambda^{n-1}|\sin\varphi|^{n-1}(1+\mu^{n+1})}\int_{-\infty}^0\frac{d\nu}{1+\nu^2}\le\frac{\operatorname{const}_\varphi}{\lambda(1+\mu^{n+1})}\,,\quad\forall\lambda\ge1.
\end{eqnarray}
In order to estimate the first integral in the RHS of (\ref{Estimand})
let us perform a change of variable $\nu\mapsto\mu+\nu$,
\begin{multline}
\label{dima111}
\int_{-\infty}^{+\infty}\left[g_{n-1}(\nu,\lambda e^{i\varphi})
-g_{n-1}(\mu,\lambda e^{i\varphi})\right]\rho(\nu-\mu)\,d\nu\\
=\int_{-\infty}^{+\infty}\left[g_{n-1}(\mu+\nu,\lambda e^{i\varphi})
-g_{n-1}(\mu,\lambda e^{i\varphi})\right]\rho(\nu)\,d\nu.
\end{multline}
Writing Taylor's formula with remainder in Lagrange's form and using \eqref{BPartDerivProp},
we get
\begin{multline}
\label{BTaylor}
g_{n-1}(\mu+\nu,\lambda e^{i\varphi})-g_{n-1}(\mu,\lambda e^{i\varphi})=
-(n-1)g_{n}(\mu,\lambda e^{i\varphi})\,\nu
\\
+\frac{n(n-1)}2g_{n+1}(\mu,\lambda e^{i\varphi})\,\nu^2
-\frac{(n+1)n(n-1)}6R(\mu,\nu,\lambda,\varphi)\,\nu^3,
\end{multline}
where
\begin{equation}
\label{formula for remainder}
R(\mu,\nu,\lambda,\varphi)
=g_{n+2}(\xi_{\mu,\mu+\nu},\lambda e^{i\varphi})
\end{equation}
and $\xi_{\mu,\mu+\nu}$ is some real number strictly between $\mu$ and $\mu+\nu$.
From (\ref{rhoProp}), (\ref{BTaylor}) and (\ref{BPartDerivProp}) we obtain
\begin{multline}
\label{dima222}
\int_{-\infty}^{+\infty}\left[g_{n-1}(\mu+\nu,\lambda e^{i\varphi})-g_{n-1}(\mu,\lambda e^{i\varphi})\right]
\rho(\nu)\,d\nu
\\
=-\frac{(n+1)n(n-1)}6\int_{-\infty}^{+\infty}R(\mu,\nu,\lambda,\varphi)\,\nu^3\rho(\nu)\,d\nu.
\end{multline}

Comparing formula \eqref{ToProveThat}
with
\eqref{Estimand}--\eqref{dima111}
and \eqref{dima222}
we see that the proof
of Lemma~\ref{Strategy lemma 1}
has been reduced to proving that
\begin{equation}
\label{dima333}
\int_{-\infty}^{+\infty}
\left|
R(\mu,\nu,\lambda,\varphi)\,\nu^3\rho(\nu)\right|\,d\nu
\le\frac{\operatorname{const}_\varphi}{\lambda(1+\mu^{n+1})}\,,
\quad\forall\lambda\ge1,\quad\forall\mu\ge0.
\end{equation}

In order to prove \eqref{dima333}
it is sufficient to prove the following two estimates:
\begin{equation}
\label{dima444}
\int_{-\infty}^{+\infty}
\left|
R(\mu,\nu,\lambda,\varphi)\,\nu^3\rho(\nu)\right|\,d\nu
\le\frac{\operatorname{const}_\varphi}{\lambda^{n+2}}\,,
\quad\forall\lambda\ge1,\quad\forall\mu\in[0,\lambda],
\end{equation}
\begin{equation}
\label{dima555}
\int_{-\infty}^{+\infty}
\left|
R(\mu,\nu,\lambda,\varphi)\,\nu^3\rho(\nu)\right|\,d\nu
\le\frac{\operatorname{const}_\varphi}{\lambda\mu^{n+1}}\,,
\quad\forall\lambda\ge1,\quad\forall\mu\ge\lambda.
\end{equation}

Observe that
formulae
\eqref{formula for remainder}
and
\eqref{BBound}
give us the rough estimate
\begin{equation}
\label{dima123}
|R(\mu,\nu,\lambda,\varphi)|
\le\frac{6}{|\sin\varphi|^{n+2}\lambda^{n+2}}\,,
\quad\forall\lambda>0,\quad\forall\mu,\nu\in\mathbb{R}.
\end{equation}
Formulae
\eqref{dima123}
and
\eqref{rhoSchwartz} with $p=5$
imply 
\eqref{dima444}.

Formulae
\eqref{formula for remainder}
and
\eqref{BBound}
also
tell us that
\[
|R(\mu,\nu,\lambda,\varphi)|
\le\frac{\operatorname{const}_\varphi}{\mu^{n+2}}
\le\frac{\operatorname{const}_\varphi}{\lambda\mu^{n+1}}
\]
uniformly over all $\mu\ge\lambda>0$ and $\nu\ge-\mu/2$.
Using this estimate and formula
\eqref{rhoSchwartz} with $p=5$
we get
\begin{equation}
\label{dima666}
\int_{-\mu/2}^{+\infty}
\left|
R(\mu,\nu,\lambda,\varphi)\,\nu^3\rho(\nu)\right|\,d\nu
\le\frac{\operatorname{const}_\varphi}{\lambda\mu^{n+1}}\,,
\quad\forall\lambda\ge1,\quad\forall\mu\ge\lambda.
\end{equation}

Comparing formulae
\eqref{dima666}
and
\eqref{dima555}
we see that the proof
of Lemma~\ref{Strategy lemma 1}
has been reduced to proving that
\begin{equation}
\label{dima777}
\int_{-\infty}^{-\mu/2}
\left|
R(\mu,\nu,\lambda,\varphi)\,\nu^3\rho(\nu)\right|\,d\nu
\le\frac{\operatorname{const}_\varphi}{\lambda\mu^{n+1}}\,,
\quad\forall\lambda\ge1,\quad\forall\mu\ge\lambda.
\end{equation}

Using
\eqref{dima123}
and
\eqref{rhoSchwartz}
with $p=n+5$ we get
\begin{multline*}
\int_{-\infty}^{-\mu/2}\,
\left|
R(\mu,\nu,\lambda,\varphi)\,\nu^3\rho(\nu)
\right|\,d\nu
\le
\frac{6c_{n+5}}{|\sin\varphi|^{n+2}\lambda^{n+2}}
\int_{\mu/2}^{+\infty}\frac{d\nu}{\nu^{n+2}}
\\
=
\frac{6\cdot2^{n+1}c_{n+5}}{(n+1)|\sin\varphi|^{n+2}\lambda^{n+2}\mu^{n+1}}
\,,
\quad\forall\lambda\ge1,\quad\forall\mu\ge\lambda,
\end{multline*}
which implies \eqref{dima777}.
\qed

\section{Some integrals involving the functions~$g_n$}
\label{Some integrals involving the functions}

In this appendix we evaluate some integrals involving the functions \eqref{definition of g}.
These results will be used later in Appendix \ref{Proof of Lemma Strategy lemma 2}.

Let us evaluate the following indefinite integral:
\begin{multline}
\label{some1}
\int\frac{\mu^nd\mu}{(\mu-z)^n}=\int\left(1+\frac{z}\nu\right)^nd\nu
=\int\left[1+\frac{nz}\nu+\sum_{k=2}^n\binom{n}{k}\frac{z^k}{\nu^k}\right]d\nu
\\
=\nu+nz\log\nu+\sum_{k=2}^n\binom{n}{k}\frac1{1-k}z^k\nu^{1-k}
\\
=\mu+nz\log(\mu-z)+\sum_{k=2}^n\binom{n}{k}\frac1{1-k}z^k(\mu-z)^{1-k}.
\end{multline}
Here  in performing intermediate calculations we used the change of variable $\nu=\mu-z$.

Similarly,
\begin{multline}
\label{some2}
\int\frac{\mu^{n-1}d\mu}{(\mu-z)^n}=\int\left(1+\frac{z}\nu\right)^{n-1}\frac{d\nu}\nu
=\int\left[1+\sum_{k=1}^{n-1}\binom{n-1}{k}\frac{z^k}{\nu^k}\right]\frac{d\nu}\nu
\\
=\log\nu-\sum_{k=1}^{n-1}\binom{n-1}{k}\frac1kz^k\nu^{-k}
=\log(\mu-z)-\sum_{k=1}^{n-1}\binom{n-1}{k}\frac1kz^k(\mu-z)^{-k}.
\end{multline}

Formulae
\eqref{definition of g},
\eqref{some1}
and
\eqref{some2}
imply
\begin{multline}
\label{some3}
\int g_n(\mu,z)\,\mu^n\,d\mu=2nz\log(\mu-z)-2nz\log(\mu-2z)
\\
+2\sum_{k=2}^n\binom{n}{k}\frac1{1-k}z^k(\mu-z)^{1-k}-\sum_{k=2}^n\binom{n}{k}\frac1{1-k}2^kz^k(\mu-2z)^{1-k}
-\mbox{c.c.},
\end{multline}
\begin{multline}
\label{some4}
\int g_n(\mu,z)\,\mu^{n-1}\,d\mu=2\log(\mu-z)-\log(\mu-2z)
\\
-2\sum_{k=1}^{n-1}\binom{n-1}{k}\frac1kz^k(\mu-z)^{-k}+\sum_{k=1}^{n-1}\binom{n-1}{k}\frac1k2^kz^k(\mu-2z)^{-k}
-\mbox{c.c.}.
\end{multline}

Using
\eqref{some3}
and
\eqref{some4}
we can finally evaluate definite integrals:
\begin{equation}
\label{some5}
\int_0^{+\infty}
g_n(\mu,z)\,\mu^n\,d\mu
=
\left.\left[
2nz\log\left(\frac{\mu-z}{\mu-2z}\right)
-
2n\bar z\log\left(\frac{\mu-\bar z}{\mu-2\bar z}\right)
\right]\right\rvert_0^{+\infty},
\end{equation}
\begin{equation}
\label{some6}
\int_0^{+\infty}
g_n(\mu,z)\,\mu^{n-1}\,d\mu
=
\left.\left[
\log\left(\frac{\mu-z}{\mu-2z}\right)
+
\log\left(\frac{\mu-z}{\mu-\bar z}\right)
-
\log\left(\frac{\mu-\bar z}{\mu-2\bar z}\right)
\right]\right\rvert_0^{+\infty}.
\end{equation}
Here the complex logarithms $\,\log\,$ are continuous multivalued functions
which have to be handled carefully.

Note that for any $z\in\mathbb{C}\setminus\mathbb{R}$
and any real positive $\mu$ we have
\[
\operatorname{Im}\frac{\mu-z}{\mu-2z}=\frac{\mu\operatorname{Im}z}{|\mu-2z|^2}\neq0,
\]
\[
\operatorname{Im}\frac{\mu-z}{\mu-\bar z}=\frac{2\operatorname{Im}z(\operatorname{Re}z-\mu)}{|\mu-z|^2}=0
\quad\Rightarrow\quad
\operatorname{Re}\frac{\mu-z}{\mu-\bar z}=\frac{(\operatorname{Re}z-\mu)^2-(\operatorname{Im}z)^2}{|\mu-z|^2}<0,
\]
so neither of the two arguments of our $\,\log\,$
crosses the positive real axis $\mathbb{R}_+\,$.
Hence, we are free to switch from $\,\log\,$ to the single-valued $\,\operatorname{Log}:\mathbb{C}\setminus\{0\}\to\mathbb{R}+i[0,2\pi)\,$ branch-cut along $\mathbb{R}_+$.
Formulae
\eqref{some5}
and
\eqref{some6}
become
\begin{multline}
\label{IntBnmu^n}
\int_0^{+\infty}
g_n(\mu,z)\,\mu^n\,d\mu
=
\left.\left[
2nz\operatorname{Log}\left(\frac{\mu-z}{\mu-2z}\right)
-
2n\bar z\operatorname{Log}\left(\frac{\mu-\bar z}{\mu-2\bar z}\right)
\right]\right\rvert_0^{+\infty}
\\
=2n(z-\bar z)\ln2
=4ni(\ln 2)\operatorname{Im}z,
\end{multline}
\begin{multline}
\label{IntBnmu^n-1}
\int_0^{+\infty} g_n(\mu,z)\,\mu^{n-1}\,d\mu
=
\left.\left[
\operatorname{Log}\left(\frac{\mu-z}{\mu-2z}\right)
+
\operatorname{Log}\left(\frac{\mu-z}{\mu-\bar z}\right)
-
\operatorname{Log}\left(\frac{\mu-\bar z}{\mu-2\bar z}
\right)\right]\right\rvert_0^{+\infty}
\\
=
\left.\operatorname{Log}\left(\frac{\mu-z}{\mu-\bar z}\right)\right\rvert_0^{+\infty}
=
i\pi(1+\operatorname{sgn}\operatorname{Im}z)-i\operatorname{Arg}z^2,
\end{multline}
where $\operatorname{Arg}:\mathbb{C}\setminus\{0\}\to[0,2\pi)$ is also branch-cut along $\mathbb{R}_+\,$.

\section{Proof of Lemma \ref{Strategy lemma 2}}
\label{Proof of Lemma Strategy lemma 2}

Formula
\eqref{expansion for mollified derivative of counting function}
tells us that
\begin{equation}
(N_\pm'*\rho)(x,\mu)=a_{n-1}^\pm(x)\mu^{n-1}+a_{n-2}^\pm(x)\mu^{n-2}+(1+\mu)^{n-3}r^\pm(x,\mu),
\end{equation}
where $r^\pm(x,\mu)$ is bounded uniformly in $\mu\ge0$.

Let $g_n(\mu,z)$ be defined in accordance with \eqref{definition of g}.
We have
\begin{equation}
\label{BHomogenProp}
g_n(\lambda\mu,\lambda z)=\lambda^{-n}g_n(\mu,z),\quad\forall\lambda>0.
\end{equation}
Using (\ref{BHomogenProp}) we get
\begin{equation}
\label{b3}
\int_0^{+\infty}g_{n-1}(\mu,\lambda e^{i\varphi})\,\mu^{n-1}\,d\mu
=\lambda\int_0^{+\infty}g_{n-1}(\mu,e^{i\varphi})\,\mu^{n-1}\,d\mu\,,
\end{equation}	
\begin{equation}
\label{b4}
\int_0^{+\infty}g_{n-1}(\mu,\lambda e^{i\varphi})\,\mu^{n-2}\,d\mu
=\int_0^{+\infty}g_{n-1}(\mu,e^{i\varphi})\,\mu^{n-2}\,d\mu,
\end{equation}	
\begin{multline}
\label{b5}
\int_0^{+\infty}g_{n-1}(\mu,\lambda e^{i\varphi})\,(1+\mu)^{n-3}\,r^\pm(x,\mu)\,d\mu
\\
=\frac1\lambda\int_0^{+\infty} g_{n-1}(\mu,e^{i\varphi})
\left(\frac1{\lambda}+\mu\right)^{n-3}
r^\pm(x,\lambda\mu)\,d\mu=o(1)
\quad
\text{as}
\quad
\lambda\to+\infty.
\end{multline}
Recall (see Appendix \ref{Proof of Lemma Strategy lemma 1})
that
the function $g_{n-1}(\mu,z)$ decays as $\mu^{-n-1}$
when $\mu\to+\infty$,
so the integrals in
\eqref{b3}--\eqref{b5}
converge.

Substituting \eqref{b3}--\eqref{b5} into (\ref{frhodef}) we get
\begin{multline}
\label{extralabel1}
f^\rho(x,\lambda e^{i\varphi})
\\
=\lambda i
\left[
a_{n-1}^+(x)\int_0^{+\infty}g_{n-1}(\mu,e^{i\varphi})\,\mu^{n-1}\,d\mu
\,-\,
(-1)^na_{n-1}^-(x)\int_0^{+\infty}g_{n-1}(\mu,e^{i(\varphi+\pi)})\,\mu^{n-1}\,d\mu
\right]
\\
+i
\left[
a_{n-2}^+(x)\int_0^{+\infty}g_{n-1}(\mu,e^{i\varphi})\,\mu^{n-2}\,d\mu
\,-\,
(-1)^na_{n-2}^-(x)\int_0^{+\infty}g_{n-1}(\mu,e^{i(\varphi+\pi)})\,\mu^{n-2}\,d\mu
\right]
\\
+o(1)
\quad
\text{as}
\quad
\lambda\to+\infty.
\end{multline}
Formulae
(\ref{IntBnmu^n})
and
(\ref{IntBnmu^n-1})
give us the values of the integrals appearing in
\eqref{extralabel1},
so
\eqref{extralabel1}
becomes
\begin{multline*}
f^\rho(x,\lambda e^{i\varphi})
=
-4(n-1)(\ln2)(\sin\varphi)\left[a_{n-1}^+(x)+(-1)^na_{n-1}^{-}(x)\right]\lambda
\\
-2\left[a_{n+2}^+(x)(\pi-\varphi)+(-1)^na_{n-1}^-(x)\varphi\right]
+o(1)
\quad
\text{as}
\quad
\lambda\to+\infty,
\end{multline*}
thus proving the lemma.
\qed

\section{Weyl quantization on manifolds}
\label{Weyl quantization on manifolds}\footnote{The
    content of this appendix can be found in a slightly more concentrated
    form in the appendix of \cite{SjZw02}.
    The main ideas and related results appeared earlier in Appendix a.3 in \cite{HeSj89}. We
  recovered these precise references only after completing the section
and decided to keep it for the convenience of the
  reader. See also Section 18.5 in \cite{Ho85}. }\label{avantpropos}
Let $M$ be a compact manifold. A pseudodifferential operator of order
$m\in \mathbb{R}$ is a continuous operator $A:C^\infty (M)\to C^\infty
(M)$ which has a weakly continuous extension ${\cal D}'(M)\to {\cal
  D}'(M)$ such that, with $K_A$ denoting the distribution kernel,
\begin{itemize}
\item[1)] $\mathrm{sing\, supp}K_A\subset \mathrm{diag\,}(M\times
  M)$,
\item[2)] For every system of local coordinates $\gamma :\Omega \ni
  \rho \mapsto 
x\in \Omega '\subset \mathbb{R}^n
$
where $\Omega \subset M$, $\Omega '\subset \mathbb{R}^n$ are open and $\gamma
$ a diffeomorphism, we have (identifying $\Omega $ and $\gamma (\Omega)$)
\begin{equation}\label{D.1}
Au(x)=\frac{1}{(2\pi )^n}\iint e^{i(x-y)\cdot \theta } a(x,\theta
)u(y)dyd\theta +Ru,\ \ u\in C_0^\infty (\Omega ),\ x\in \Omega ,
\end{equation}
where $R$ is smoothing ($K_R\in C^\infty (\Omega \times \Omega ) $) and
$a$ is a symbol of order $m$; $a\in S^m(\Omega)$, which means that
$a\in C^\infty (\Omega \times \mathbb{R}^n)$ and that for every $\widehat{K}\Subset
\Omega $ and all $\alpha ,\beta \in \mathbb{N}^n$, $\exists C=C_{\widehat{K},\alpha
,\beta }$ such that
\begin{equation}\label{D.2}
|\partial _x^\alpha \partial _\theta ^\beta a(x,\theta )|\le C\langle
\theta \rangle^{m-|\beta |},\ \forall\, (x,\theta )\in
\widehat{K}\times \mathbb{R}^n , \hbox{ where }\langle \theta \rangle=
(1+|\theta |^2)^{1/2}.
\end{equation}
\end{itemize}  

\par If $\widetilde{\gamma }:\widetilde{\Omega }\ni \rho \mapsto
\widetilde{x}\in \widetilde{\Omega }'$ is
another local coordinate chart, then over the intersection $\Omega
\cap \widetilde{\Omega }$ we can express $x=\kappa
(\widetilde{x})$, where $\kappa =\gamma \circ
\widetilde{\gamma }^{-1}$ and we have
\begin{equation}\label{D.3}
a(\kappa (\widetilde{x}),\theta )\equiv
\widetilde{a}(\widetilde{x},(\kappa '(\widetilde{x}))^\mathrm{t}\theta
)\ \mathrm{mod\,}S^{m-1}.
\end{equation}
This allows us to define the symbol $\sigma _A$ of $A$ on $T^*M$ up to
symbols of order 1 lower. More precisely, we have a bijection
\begin{equation}\label{D.4}
L^m(M)/L^{m-1}(M)\ni A\mapsto \sigma _A\in S^m(T^*M)/S^{m-1}(T^*M),
\end{equation}
with the natural definition of the symbol classes $S^m(T^*M)$, and with
$L^m(M)$ denoting the space of pseudodifferential operators on $M$ of order
$m$.

\par It is well known that we can 
replace $a(x,\theta )$ in (\ref{D.1}) with $a((x+y)/2,\theta )$ and
this leads to the same definition of $\sigma _A$ in
$S^m/S^{m-1}(T^*M)$. Thus, working with 
\begin{equation}\label{D.5}
Au(x)=\mathrm{Op\,}(a)u(x)+Ru, \ \ a\in S^m(\Omega \times {\bf
  R}^{n}),\ K_R\in C^\infty ,
\end{equation}
leads to the same principal symbol map. Here we
write\footnote{Strictly speaking, when $\Omega $ is not convex we
  need here to insert a suitable smooth cutoff $\chi (x,y)\in C^\infty
  (\Omega \times \Omega )$ which is equal to one near the diagonal,
  the choice of which can affect the operator only by a smoothing one.}
\begin{equation}\label{D.6}
\mathrm{Op\,}(a)u(x)=\frac{1}{(2\pi )^n}\iint e^{i(x-y)\cdot \theta
}a\left(\frac{x+y}{2},\theta  \right) u(y)dyd\theta .
\end{equation}

It seems to be a well-known result (though we did not find a precise
reference) that if we fix a positive smooth density $\omega $ on $M$,
restrict our attention to local coordinates for which $\omega
=dx_1...dx_n$ and work with the Weyl quantization as in (\ref{D.5}),
(\ref{D.6}), then (\ref{D.4}) improves to a bijection
\begin{equation}\label{D.7}
L^m/L^{m-2}(M)\ni A\mapsto \sigma _A\in S^m/S^{m-2}(T^*M).
\end{equation}

A natural generalization of this is to consider pseudodifferential
operators acting on $1/2$-densities; $A:C^\infty (M;\Omega ^{1/2})\to
C^\infty (M;\Omega ^{1/2}) $. When using the Weyl quantization, we get
the local representation analogous to \ref{D.5}:
\begin{equation}\label{D.8}
A(u(y)dy^{1/2})=(\mathrm{Op\,}(a)u)(x)dx^{1/2}+(Ru)dx^{1/2},
\end{equation}
where $dx=dx_1...dx_n$. Recall that Duistermaat and H\"ormander
\cite{DuiHor} have defined invariantly the notion of subpincipal symbol
of such operators when the symbols are sums of a leading
positively homogeneous term of order $m$ in $\xi $ and a symbol of
order $m-1$. This result, as well as the fixed density invariance
mentioned above, follow from the 
next more or less well-known proposition (cf.~the footnote on page \pageref{avantpropos}).
\begin{proposition}\label{D1}
Let $L^m(M)$ denote the space of pseudodifferential operators on $M$
of order $m$, acting on half densities. Then if $(x_1,...,x_n)$ and
$(\widetilde{x}_1,...,\widetilde{x}_n)$ are two local coordinate
charts and we use the representation (\ref{D.8}), so that 
$$
A(udx^{1/2})\equiv (\mathrm{Op\,}(a)u)dx^{1/2}\equiv
(\mathrm{Op\,}(\widetilde{a})\widetilde{u})d\widetilde{x}^{1/2}, 
$$
modulo the action of smoothing operators, for
$udx^{1/2}=\widetilde{u}d\widetilde{x}^{1/2}$ supported in the
intersection of the two coordinate charts, then we have
\begin{equation}\label{D.9}
a(\kappa (\widetilde{x}),\theta ) \equiv
\widetilde{a}(\widetilde{x},\kappa '(\widetilde{x})^\mathrm{t}\theta
)\ \mathrm{mod}\ S^{m-2},
\end{equation}
implying that we have a natural bijective symbol map
\begin{equation}\label{D.10}
L^m/L^{m-2}(M)\to S^m/S^{m-2}(T^*M).
\end{equation}
\end{proposition}

\par\noindent {\bf Proof. } We only verify (\ref{D.9}) and omit the
(even more) standard arguments for (\ref{D.10}). Our proof will be a
straightforward adaptation of the proof of the invariance of
pseudo\-differential operators under composition with diffeomorphisms by
means of the Kuranishi trick (cf.\ \cite{GrSj94}).

\par In the intersection of the two coordinate charts $\Omega $ and
$\widetilde{\Omega }$, we have
$u(y)dy^{1/2}=\widetilde{u}(\widetilde{y})d\widetilde{y}^{1/2}$. Here
$y=\kappa (\widetilde{y})$, where $\kappa $ is a diffeomorphism:
$\widetilde{\gamma }(\Omega \cap \widetilde{\Omega })\to \gamma
(\Omega \cap \widetilde{\Omega })$, $\kappa =\gamma \circ
\widetilde{\gamma }^{-1}$). Thus
$u(y)=\widetilde{u}(\widetilde{y})(\det \kappa
'(\widetilde{y}))^{-1/2}$, assuming that $\det \kappa '>0$ for
simplicity. Thus, modulo the action of smoothing operators
$$
A(udy^{1/2})\equiv (\mathrm{Op\,}(a)u)dx^{1/2}=(\det \kappa '(\widetilde{x}))^{1/2}(\mathrm{Op\,}(a)u)d\widetilde{x}^{1/2},
$$
so up to a smoothing operator $\mathrm{Op\,}(\widetilde{a})$ coincides
with
$$
B:\, \widetilde{u}\mapsto (\det \kappa
'(\widetilde{x}))^{1/2}\mathrm{Op\,}(a)u,\ \
u(y)=\widetilde{u}(\widetilde{y})(\det \kappa '(\widetilde{y}))^{-1/2}.
$$
We have
\begin{multline*}
B\widetilde{u}(\widetilde{x})=(\det \kappa
'(\widetilde{x}))^{1/2}\iint e^{i(x-y)\cdot \theta
}a\left(\frac{x+y}{2},\theta  \right) u(y)dy\frac{d\theta}{(2\pi )^n}\\ 
=(\det \kappa '(\widetilde{x}))^{1/2}\iint e^{i(\kappa
  (\widetilde{x})-y)\cdot \theta }a\left(\frac{\kappa (\widetilde{x})+y}{2},\theta
\right) \widetilde{u}(\widetilde{y})(\det \kappa
'(\widetilde{y}))^{-1/2}dy\frac{d\theta }{(2\pi )^n}\\
=\iint e^{i(\kappa (\widetilde{x})-\kappa (\widetilde{y}))\cdot \theta
}a\left(\frac{\kappa (\widetilde{x})+\kappa (\widetilde{y})}{2},\theta
  \right)\widetilde{u}(\widetilde{y})(\det \kappa '(\widetilde{x})
  \det \kappa '(\widetilde{y}))^{1/2}
  d\widetilde{y}\frac{d\theta}{(2\pi )^n}\,.
\end{multline*}

\par By Taylor's formula (and restricting to a suitably thin
neighborhood of the diagonal by means of a smooth cutoff, equal to one
near the diagonal), we get
$$
\kappa (\widetilde{x})-\kappa (\widetilde{y})=K(\widetilde{x},\widetilde{y})(\widetilde{x}-\widetilde{y}),
$$ 
where $\widetilde{K}(\widetilde{x},\widetilde{y})$ depends smoothly on
$(\widetilde{x},\widetilde{y})$ and 
$$
K(\widetilde{x},\widetilde{y})=\kappa '
\left(\frac{\widetilde{x}+\widetilde{y}}{2} \right) +{O}((\widetilde{x}-\widetilde{y})^2).
$$
It follows that
\begin{multline*}
B\widetilde{u}(\widetilde{x})=\iint
e^{i(\widetilde{x}-\widetilde{y})\cdot
  K^{\mathrm{t}}(\widetilde{x},\widetilde{y})\theta
}a\left(\frac{\kappa (\widetilde{x})+\kappa (\widetilde{y})}{2},\theta
\right) \widetilde{u}(\widetilde{y})(\det \kappa '(\widetilde{x}) \det
\kappa '(\widetilde{y}))^{1/2}d\widetilde{y}\frac{d\theta}{(2\pi
  )^n}\\
= \iint
e^{i(\widetilde{x}-\widetilde{y})\cdot
  \widetilde{\theta }
}a\left(\frac{\kappa (\widetilde{x})+\kappa
    (\widetilde{y})}{2},K^{\mathrm{t}}(\widetilde{x},\widetilde{y})^{-1}\widetilde{\theta }
\right) \widetilde{u}(\widetilde{y})\frac{(\det \kappa '(\widetilde{x}) \det
\kappa '(\widetilde{y}))^{1/2}}{\det K(\widetilde{x},\widetilde{y})} d\widetilde{y}\frac{d\widetilde{\theta }}{(2\pi
  )^n}.
\end{multline*}
Here
\[
\begin{split}
\frac{\kappa (\widetilde{x})+\kappa (\widetilde{y})}{2}&=\kappa
\left(\frac{\widetilde{x}+\widetilde{y}}{2} \right)
+{O}((\widetilde{x}-\widetilde{y})^2),\\
K^{\mathrm{t}}(\widetilde{x},\widetilde{y})^{-1}&=\left((\kappa ')^{\mathrm{t}}\left(\frac{\widetilde{x}+\widetilde{y}}{2} \right)  \right)^{-1}
+{O}((\widetilde{x}-\widetilde{y})^2),\\
\det K(\widetilde{x},\widetilde{y})&=\det \kappa '\left(\frac{\widetilde{x}+\widetilde{y}}{2} \right)
+{O}((\widetilde{x}-\widetilde{y})^2),\\
(\det \kappa '(\widetilde{x}) \det
\kappa '(\widetilde{y}))^{1/2}
&=\det \kappa '\left(\frac{\widetilde{x}+\widetilde{y}}{2} \right)
+{O}((\widetilde{x}-\widetilde{y})^2).
\end{split}
\]
Thus,
$$
B\widetilde{u}=\mathrm{Op\,}(\widetilde{a})\widetilde{u}+\iint
e^{i(\widetilde{x}-\widetilde{y})\cdot \widetilde{\theta
  }}b(\widetilde{x},\widetilde{y},\widetilde{\theta
})u(\widetilde{y})d\widetilde{y}\frac{d\widetilde{\theta }}{(2\pi )^n},
$$
where $\widetilde{a}\in S^m$ is related to $a$ as in (\ref{D.9}) and
$b\in S^m(\widetilde{\gamma }(\Omega \cap \widetilde{\Omega
})^2\times \mathbb{R}^n)$ (in the sense that $\partial
_{\widetilde{x}}^\alpha \partial _{\widetilde{y}}^\beta \partial
_{\widetilde{\theta }}^{|\delta|} b={O}(\langle \widetilde{\theta }
\rangle^{m-\delta })$ uniformly in $\widetilde{\theta }$ and locally
uniformly in $(\widetilde{x},\widetilde{y})$) and $b$ vanishes to the
second order on the diagonal, $\widetilde{x}=\widetilde{y}$. By
standard arguments we have $B\equiv \mathrm{Op}(r)$, where $r\in
S^{m-2}$ and the proposition follows. \hfill{$\Box$}

\section{The resolvent and its powers as
    pseudo\-differential operators}
\label{Symbolic approximation for the resolvent and its powers}

Let $\gamma :M\supset \Omega \to \Omega '\subset \mathbb{R}^n$ be a chart
of local coordinates and let us identify $\Omega '$ with $\Omega $ in the
natural way. Let $a(x,\xi )\in S^1(\Omega \times \mathbb{R}^n)$ (defined
modulo $S^{-\infty }(\Omega \times \mathbb{R}^n)$) be the Weyl symbol of
\begin{equation}\label{E.1}{{A}_\vert}_{C_0^\infty (\Omega )}:C_0^\infty (\Omega )\to C^\infty
(\Omega ), \end{equation} so that 
\begin{equation}\label{E.2}
Au(x)=\mathrm{Op\,}(a)u(x)+Ru(x),\ x\in \Omega 
\end{equation}
for every $u\in C_0^\infty (\Omega )$, where $R\in L^{-\infty }(\Omega
)$ in the sense that $K_R\in C^\infty (\Omega \times \Omega )$. Here
we identify 1/2 densities and scalar functions on $\Omega $ by means
of the fixed factor $dx^{1/2}$. We
first work in this fixed local coordinate chart and write simply $A$
for the operator in (\ref{E.1}). We notice that 
\begin{equation}\label{E.3}
a-z\in S(\Omega \times \mathbb{R}^n,\langle \xi ,z\rangle )=S(\langle \xi
,z\rangle),
\end{equation} 
in the sense that $a-z\in C^\infty (\Omega \times \mathbb{R}^n)$ and that
for all $K\Subset \Omega $, $\alpha ,\beta \in \mathbb{N}^n$,
\begin{equation}\label{E.4}
|\partial _x^\alpha \partial _\xi ^\beta (a-z)|\le C_{K,\alpha ,\beta
}\langle \xi ,z\rangle \langle \xi \rangle^{-|\beta |},
\end{equation}
uniformly when $z\in \mathbb{C}$, $|z|>1$, $x\in K$, $\xi\in \mathbb{R}^n
$. Here, we write $\langle \xi \rangle = (1+|\xi |^2)^{1/2}$, $\langle
\xi ,z\rangle =(1+|z|^2+|\xi |^2)^{1/2}$. 

\par Similarly, if $\Gamma \subset \dot{\mathbb{C}}$ is a closed conic
neighborhood of $\dot{\mathbb{R}}$ and until further notice we restrict our attention to
$z\in\dot{\mathbb{C}}\setminus (\Gamma \cup D(0,1)) $, we have 
\begin{equation}\label{E.5}
(a-z)^{-1}\in S(\langle \xi ,z\rangle^{-1})
\end{equation}
with the natural generalization of the definition (\ref{E.4}). 

Sometimes, we shall exploit the fact that $a-z$ and $(a-z)^{-1}$ belong to
narrower symbol classes, used in \cite{GrSe95}. We say that $b(x,\xi,
z) $, defined for $(x,\xi ,z)$ as in (\ref{E.5}), belongs to
$S_1(\langle \xi ,z\rangle ^m)$, $m\in \mathbb{R}$, if 
\begin{equation}\label{E.6}
|\partial _x^\alpha \partial _\xi ^\beta b(x,\xi ,z)|\le C_{K,\alpha
  ,\beta }\begin{cases}
\langle \xi ,z\rangle ^m, \hbox{ when }\alpha =\beta =0,\\
\langle \xi ,z\rangle ^m\frac{\langle \xi \rangle}{\langle \xi
  ,z\rangle}\langle \xi \rangle^{-|\beta |},\hbox{ when }(\alpha
,\beta )\ne (0,0),
\end{cases}
\end{equation}
uniformly for $x\in K\Subset \Omega $, $\xi \in \mathbb{R}^n$, $z\in
\dot{\mathbb{C}}\setminus (\Gamma \cup D(0,1))$.

If $b_j\in S(\langle \xi ,z\rangle^{m_j})$, $j=1,2$, the asymptotic
Weyl composition
\begin{multline}\label{E.7}
b_1\# b_2=\left(e^{(i/2)\sigma (D_{x,\xi };D_{y,\eta })}b_1(x,\xi
)b_2(y,\eta )\right)_{y=x\atop \eta =\xi }\\
\sim \sum_{k=0}^\infty \frac{1}{k!}\left(\left(\frac{i}{2}\sigma (D_{x,\xi
  };D_{y,\eta }) \right)^kb_1(x,\xi )b_2(y,\eta )\right)_{y=x\atop \eta =\xi }
\end{multline}
is well defined in $S(\langle \xi ,z\rangle^{m_1+m_2})/S(\langle \xi
,z\rangle^{m_1+m_2}\langle \xi \rangle^{-\infty })$, where $$S(\langle \xi
,z\rangle^{m_1+m_2}\langle \xi \rangle^{-\infty })=\bigcap_{N\ge 0}S(\langle \xi
,z\rangle^{m_1+m_2}\langle \xi \rangle^{-N})$$ and with the natural
definition of the symbol spaces to the right.
Here $\sigma (D_{x,\xi};D_{y,\eta })=D_\xi \cdot D_y-D_x\cdot D_\eta $. Notice that
\begin{equation}\label{E.8}
\frac{1}{k!}\left(\left(\frac{i}{2}\sigma (D_{x,\xi
  };D_{y,\eta }) \right)^kb_1(x,\xi )b_2(y,\eta )\right)_{y=x\atop
\eta =\xi }
\in S(\langle \xi ,z\rangle^{m_1+m_2}\langle \xi \rangle^{-k}).
\end{equation}

\par When $b_j\in S_1(m_j)$ this improves to
\begin{equation}\label{E.9}
\frac{1}{k!}\left(\left(\frac{i}{2}\sigma (D_{x,\xi
  };D_{y,\eta }) \right)^kb_1(x,\xi )b_2(y,\eta )\right)_{y=x\atop
\eta =\xi }
\in\begin{cases}
S_1(\langle \xi ,z\rangle^{m_1+m_2}),\hbox{ when }k=0,\\
S(\langle \xi ,z\rangle^{m_1+m_2-2}\langle \xi \rangle^{2-k}),\hbox{
  when }k\ge 1.
\end{cases}
\end{equation}
In particular,
$$
b_1\# b_2\equiv b_1b_2\hbox{ mod }S(\langle \xi
,z\rangle^{m_1+m_2-2}\langle \xi \rangle ).
$$
In the special case $b_1=a-z$, $b_2=(a-z)^{-1}$ we get
\begin{equation}\label{E.10}
(a-z)\# (a-z)^{-1}=1+r,
\end{equation}
\begin{equation}\label{E.11}\begin{split}
      r\sim & \sum_{k=1}^\infty \frac{1}{k!}\left(\frac{i}{2}\sigma (D_{x,\xi
        };D_{y,\eta }) \right)^k\left( a(x,\xi )(a(y,\eta )-z)^{-1}
      \right)_{y=x\atop \eta =\xi }\\ & \in S(\langle \xi ,z\rangle^{-2}\langle
      \xi \rangle)/S(\langle \xi ,z\rangle^{-2}\langle
      \xi \rangle^{-\infty }),\\
      r\equiv & \frac{i}{2}\sigma (D_{x,\xi };D_{y,\eta })\left(a(x,\xi)(a(y,\eta )-z)^{-1} \right)_{y=x\atop \eta =\xi }\ \\
      \equiv &\frac{i}{2}\{ a,(a-z)^{-1} \}\ \mathrm{mod}\
      S(\langle \xi ,z\rangle^{-2}),
\end{split}
\end{equation}
with the Poisson bracket as defined in Section \ref{Statement of the problem}.

\par The symbolic inverse of $A-z$ is now
\begin{equation}\label{E.12}
b(x,\xi ,z)\sim (a-z)^{-1}\# (1-r+r\#r...\pm r^{\# k}+...),
\end{equation}
where
$$r^{\#k}=\underbrace{r\#r\#...\#r}_{k\ \mathrm{ factors}}\in S\left((\langle \xi \rangle/\langle \xi
,z\rangle^2)^k\right)\subset S\left(\langle \xi ,z\rangle ^{-k}\right). $$
We see that $b(x,\xi ,z)\in S(\langle \xi ,z\rangle^{-1})$ and that 
$$
b\equiv (a-z)^{-1}\ \mathrm{mod}\ S\left(\frac{\langle \xi
    \rangle}{\langle \xi ,z\rangle^3} \right).
$$
More precisely,
$$
b\equiv (a-z)^{-1}-(a-z)^{-1}\# r\ \mathrm{mod}\
S\left(\frac{1}{\langle \xi ,z\rangle^3} \right) .
$$
Here
\begin{multline*}
(a-z)^{-1}\# r\sim (a-z)^{-1}r+\sum_{k\ge
  1}\frac{1}{k!}\left(\left(\frac{i}{2}\sigma (D_{x,\xi };D_{y,\eta })
  \right)^k\left((a-z)^{-1}(x,\xi )r(y,\eta ) \right)
\right)_{y=x\atop \eta =\xi }\\
\equiv (a-z)^{-1}r\ \mathrm{mod}\ S\left( \frac{1}{\langle \xi ,z\rangle^3} \right),
\end{multline*}
so
\begin{equation}\label{E.13}\begin{split}
b(x,\xi ,z)&\equiv (a-z)^{-1}-(a-z)^{-1}r \\
&\equiv (a-z)^{-1}-\frac{i}{2}(a-z)^{-1}\{ a, (a-z)^{-1} \}\ \mathrm{mod}\ S\left(
  \frac{1}{\langle \xi ,z\rangle^3} \right), 
\end{split}
\end{equation}
where we also used (\ref{E.11}).

\par If $b_j\in S(\langle \xi ,z\rangle^m\langle \xi \rangle ^{k-j})$
for $j=0,1,...$, we can apply a standard procedure to construct a
symbol $b\in S((\langle \xi ,z\rangle ^m\langle \xi \rangle^k)$ such
that 
$$
b-\sum_0^{N-1}b_j\in S(\langle \xi ,z\rangle^m \langle \xi \rangle^{k-N})
$$
for every $N\ge 1$ and we still write $b\sim \sum_0^\infty b_j$ where
$b$ is a concrete symbol (uniquely determined up to $S(\langle \xi
,z\rangle^m\langle \xi \rangle^{-\infty })$). If $b_j$ are holomorphic
for $z\in \dot{\mathbb{C}}\setminus (\Gamma \cup D(0,1))$, then the standard
construction produces a symbol $b$ which is also holomorphic. 

\par If $b\in S(\langle \xi ,z\rangle^m\langle \xi \rangle^k)$ is such
a holomorphic symbol then by the Cauchy inequalities we
get\footnote{After replacing $\Gamma $ with any closed conic set
  containing $\Gamma $ in its interior and $D(0,1)$ with
  $D(0,1+\epsilon )$ for any $\epsilon >0$}
$$
\partial _z^\ell b\in S(\langle \xi ,z\rangle^m\langle \xi
\rangle^k\langle z\rangle^{-\ell})
$$
in the sense that 
$$
|\partial _x^{\alpha }\partial _\xi ^{\beta }\partial _z^\ell b|
\le
C_{K,\alpha ,\beta ,\ell}\langle \xi ,z\rangle^m\langle \xi
\rangle^{k-|\beta |}\langle z\rangle^{-\ell}
$$
for $x\in K\Subset \Omega $, $\xi \in \mathbb{R}^n$ and omitting the
slight increase of $\Gamma \cup D(0,1)$, mentioned in the
last footnote.

With the holomorphic $z$-dependence in mind we return to (\ref{E.11}) and write
\begin{equation}\label{E.14}
r\sim \sum_{k=1}^\infty r_k(x,\xi ,z)
\end{equation} 
and get a concrete symbol $r\in S(\langle \xi ,z\rangle^{-2}\langle
\xi \rangle^{2-1})$ which is holomorphic in $z$, so that for every
$N\ge 1$,
\begin{equation}\label{E.15}
r-\sum_1^{N-1}r_k\in S(\langle \xi ,z\rangle^{-2}\langle \xi \rangle^{2-N})
\end{equation}
and by the Cauchy inequalities
\begin{equation}\label{E.16}
\partial _z^\ell \left(r-\sum_1^{N-1}r_k \right) \in S(\langle \xi
,z\rangle^{-2}\langle \xi \rangle^{2-N}\langle z\rangle^{-\ell}).
\end{equation}

\par From the explicit expression of the $r_k$ (or from observing
that they are defined for $z$ in $(\dot{\mathbb{C}}\setminus \Gamma )\cup
D(0,\langle \xi \rangle/C)$ when $\xi $ is large), we see that 
\begin{equation}\label{E.17}
\partial _z^\ell r_k\in S(\langle \xi ,z\rangle^{-2-\ell}\langle \xi \rangle^{2-k}),
\end{equation}
\begin{equation}\label{E.18}
\partial _z^\ell\left(\sum_1^{N-1} r_k  \right) \in S(\langle \xi
,z\rangle^{-2-\ell}\langle \xi \rangle^{2-1}).
\end{equation}
Choosing $N=\ell +1$ in (\ref{E.16}), (\ref{E.18}), we get 
\begin{equation}\label{E.19}
\partial _z^\ell r\in S(\langle \xi ,z\rangle^{-2-\ell}\langle \xi \rangle^{1}).
\end{equation}

This argument shows that (\ref{E.14}) is valid in the symbol space
$\widetilde{S}(\langle \xi ,z\rangle^{-2}\langle \xi \rangle^{2-1})$,
where we say that $c\in \widetilde{S}(\langle \xi ,z\rangle^{m}\langle
\xi \rangle^k)$ if $c(x,\xi ,z)$ is a smooth, holomorphic in $z$ and 
$$
\partial _z^\ell c\in \widetilde{S}(\langle \xi
,z\rangle^{m-\ell}\langle \xi \rangle ^k),\hbox{ for all }\ell\ge 0.
$$

\par In (\ref{E.12}) we can choose $r^{\# k}$ and the asymptotic sums
so that $b\in \widetilde{S}(\langle \xi,z\rangle^{-1})$ and so that (\ref{E.13}) improves to 
\begin{equation}\label{E.20}\begin{split}
    b(x,\xi ,z)&\equiv (a-z)^{-1}-(a-z)^{-1}r \\
    &\equiv (a-z)^{-1}-\frac{i}{2}(a-z)^{-1}\{ a, (a-z)^{-1} \}\
    \mathrm{mod}\ \widetilde{S}\left( \frac{1}{\langle \xi
        ,z\rangle^3} \right),
\end{split}
\end{equation}
where
\begin{equation}\label{E.21}
r\in \widetilde{S}(\langle \xi ,z\rangle^{-2}\langle \xi \rangle),\
(a-z)^{-1}\in \widetilde{S}(\langle \xi ,z\rangle^{-1}).
\end{equation}

\par In the main text we have 
$$
A=A_1+A_0+A_{-1},\ A_0=A_{\mathrm{sub}},
$$
where $A_j\in S(\langle \xi \rangle ^j)$ and $A_1$, $A_0$ are
positively homogeneous in $\xi $ of degree 1 and 0 respectively, in
the region $|\xi |\ge 1$. From the resolvent identity
\begin{multline*}
(a-z)^{-1}=(A_1-z)^{-1}-(A_1-z)^{-1}(a-A_1)(A_1-z)^{-1}\\
+(A_1-z)^{-1}(a-A_1)(a-z)^{-1}(a-A_1) (A_1-z)^{-1}
\end{multline*}
we infer that
$$
(a-z)^{-1}\equiv (A_1-z)^{-1}-(A_1-z)^{-1}(A_0+A_{-1})(A_1-z)^{-1}\
\mathrm{mod}\ \widetilde{S}(\langle \xi ,z\rangle^{-3}),
$$
hence,
$$
(a-z)^{-1}\equiv (A_1-z)^{-1}-(A_1-z)^{-1}A_0(A_1-z)^{-1}\
\mathrm{mod}\ \widetilde{S}(\langle \xi \rangle^{-1}\langle \xi ,z\rangle^{-2}).
$$
In particular,
$$
(a-z)^{-1}\equiv (A_1-z)^{-1}\ \mathrm{mod}\ \widetilde{S}(\langle \xi
,z\rangle^{-2})
$$
and from (\ref{E.20}) we get
\begin{multline}\label{E.21.5}
b\equiv
(A_1-z)^{-1}-(A_1-z)^{-1}A_0(A_1-z)^{-1}\\-\frac{i}{2}(A_1-z)^{-1}\{
A_1, (A_1-z)^{-1} \}\ \mathrm{mod}\
\widetilde{S}\left(\frac{1}{\langle \xi \rangle\langle \xi
    ,z\rangle^2} \right) ,
\end{multline}
which implies (\ref{12 December 2016 equation 1}) (cf.\ (\ref{12 December 2016 equation 2})).

\par By construction, $b$ is a realization of the symbolic inverse of
$a-z$:
$$
(a-z)\# b\equiv 1\ \mathrm{mod}\ \widetilde{S}(\langle \xi
,z\rangle^{-2}\langle \xi \rangle^{-\infty }). 
$$
Let $B=\mathrm{Op\,}(b):\, C_0^\infty (\Omega )\to C^\infty (\Omega )$
(where we also insert a suitable cutoff  $\in C^\infty (\Omega
\times \Omega )$, equal to 1 near $\mathrm{diag\,}(\Omega \times
\Omega )$). Then
\begin{equation}\label{E.22}
\partial _z^kB(z)={O}(\langle z\rangle^{-k_1}):\,
H_\mathrm{comp}^s(\Omega )\to H_{\mathrm{loc}}^{s+k_2}(\Omega ) \hbox{
  uniformly for }z\in \dot{\mathbb{C}}\setminus (\Gamma \cup D(0,1)),
\end{equation}
when $1+k=k_1+k_2$, $k_j\ge 0$, $s\in \mathbb{R}$.

Let $\chi ,\Phi \in C_0^\infty (\Omega )$, with $\Phi =1$
near $\mathrm{supp\,}(\chi )$. Then,
\begin{equation}\label{E.22.5}
(A-z)\Phi B\chi =\chi +R,
\end{equation}
where $R=R(z)$ is a smoothing operator: ${\cal D}'(\Omega )\to
C^\infty (\Omega )$, depending holomorphically on $z$,
such that $Ru=0$ when $\mathrm{supp\,}(u)\cap \mathrm{supp\,}(\chi
)=\emptyset $ and 
\begin{equation}\label{E.23}
\partial _z^kR={O}(\langle z\rangle^{-2-k}): H^{-s}(\Omega )\to
H_{\mathrm{loc}}^s(\Omega ),\ z\in \dot{\mathbb{C}}\setminus (\Gamma \cup D(0,1)),
\end{equation}
for all $s\in \mathbb{R}$, $k\ge 0$. We omit the standard proof of this,
based on the symbolic results above, starting with the identity
$$
(A-z)\Phi B\chi =[A,\Phi ]B\chi +\Phi (A-z)B\chi .
$$

\par Let $M\subset \bigcup_1^N\Omega _j$ be a finite covering of $M$
with coordinate charts as above. Recall that $A$ is a globally defined
pseudodifferential operator acting on 1/2 densities so we can now view
$A-z$ as acting: $C_0^\infty (\Omega _j;\Omega ^{1/2})\to C^\infty
(M;\Omega ^{1/2})$ for each $j$. We have a corresponding operator
$B_j$ (as ``$B$'' above), now acting on 1/2-densities, so that 
\begin{equation}\label{E.24}
B_j(udx^{1/2})=(\mathrm{Op\,}(b_j)u)dx^{1/2},\ u\in C_0^\infty (\Omega _j)
\end{equation} 
where $dx^{1/2}$ is the canonical (and $j$-dependent) 1/2-density on
$\Omega _j$. Let $\chi _j\in C_0^\infty (\Omega _j)$ form a partition
of unity on $M$. (\ref{E.22.5}) becomes
\begin{equation}\label{E.25}
(A-z)\Phi _jB_j\chi _j=\chi _j+R_j(z),
\end{equation}
where $R_j$ has the properties of ``$R$'' in (\ref{E.22}),
(\ref{E.23}) except for the fact that $R_j$ acts on 1/2-densities and
that we can actually define $R_j$ as an operator on $M$ such that
\begin{equation}\label{E.26}
\| \partial _z^kR_j\|_{{\cal L}(H^{-s},H^s(M))}\le C_s\langle
z\rangle^{-2-k},\ z\in \dot{\mathbb{C}}\setminus (\Gamma \cup D(0,1)).
\end{equation}
Here $H^s(M)$ denotes the Sobolev space of 1/2-densities of order
$s\in \mathbb{R}$. 

\par Let 
\begin{equation}\label{E.27}
B:=\sum \Phi _jB_j\chi _j:\ C^\infty (M;\Omega ^{1/2})\to C^\infty (M;\Omega ^{1/2}).
\end{equation}
Then 
\begin{equation}\label{E.28}
(A-z)B(z)=1+R(z),
\end{equation}
\begin{equation}\label{E.29}
R(z)=\sum R_j(z),
\end{equation}
\begin{equation}\label{E.30}
\partial _z^kB(z)={O}(\langle z\rangle^{-k_1}):\ H^s\to
H^{s+k_2},\hbox{ when }k+1=k_1+k_2,\ k_j\ge 0,
\end{equation}
\begin{equation}\label{E.31}\partial _z^k R(z)={O}_s(\langle
  z\rangle^{-2-k}):\, H^{-s}\to H^s,
\end{equation}
for all $s\in \mathbb{R}$.

On the other hand, by direct arguments, we know that $(A-z)^{-1}$ also
enjoys the properties (\ref{E.30}). Applying this operator to the left
in (\ref{E.28}), we get
\begin{equation}\label{E.32}
(A-z)^{-1}=B(z)-K(z),\ K(z)=(A-z)^{-1}R(z).
\end{equation}
Clearly, $K(z)$ also satisfies (\ref{E.31}).

\par Using the operator identity
(\ref{The matrix trace of a power of the resolvent equation 1}) in (\ref{E.32}),
we get
\begin{equation}\label{E.33}
(A-z)^{1-n}=B^{(n)}(z)-K^{(n)}(z),
\end{equation}
\begin{equation}\label{E.34}
B^{(n)}=\frac{1}{(n-2)!}\partial _z^{n-2}B(z),
\end{equation}
\begin{equation}\label{E.35}
K^{(n)}=\frac{1}{(n-2)!}\partial _z^{n-2}K(z)={O}_s(\langle
z\rangle^{-n}):\, H^{-s}\to H^s.
\end{equation}
From the last estimate it follows that $K^{(n)}$ is of trace class
with a continuous distribution kernel which is uniformly $={O}(\langle z\rangle^{-n})$.

\par Let $x_0$ be a point in a coordinate chart $\Omega =\Omega _j$
and assume for simplicity that $\chi =\chi _j$ is equal to 1 near that
point. Then near $(x_0,x_0)$ the distribution kernel of $B$
(identified locally with an operator acting on scalar functions)
coincides with that of $\mathrm{Op\,}(b)$, where $b$ satisfies
(\ref{E.20}). Consequently,
\begin{equation}\label{E.36}
B^{(n)}=\mathrm{Op\,}(b^{(n)}),
\end{equation}
\begin{equation}\label{E.37}
b^{(n)}\equiv (a-z)^{-n}-\frac{1}{(n-2)!}\partial _z^{n-2}\left(
  (a-z)^{-1}r \right)
\ \mathrm{mod}\ \widetilde{S}\left(\langle \xi ,z\rangle^{-n-1}\right).
\end{equation}

\section{Proof of formulae
\eqref{20 December 2016 equation 1}
and
\eqref{20 December 2016 equation 2}
}
\label{Proof of formulae (4.4) and (4.5)}

Formula
\eqref{basic identity for eigenprojections}
implies
\begin{equation}
\label{dPP}
(\partial P^{(k)})P^{(j)}+P^{(k)}\partial P^{(j)}
=\delta^{kj}\partial P^{(k)},
\end{equation}
where $\partial$ is any partial derivative. We have
\begin{multline*}
{\operatorname{tr}}
\{P^{(k)},P^{(j)},P^{(l)}\}
={\operatorname{tr}}
\bigl[(\partial_{x^\alpha}P^{(k)})P^{(j)}\partial_{\xi_\alpha}P^{(l)}
-
(\partial_{\xi_\alpha}P^{(k)})P^{(j)}\partial_{x^\alpha}P^{(l)}\bigr]
\\
=
{\operatorname{tr}}
\bigl[
\bigl((\partial_{x^\alpha}P^{(k)})P^{(j)}\bigr)\bigl(P^{(j)}\partial_{\xi_\alpha}P^{(l)}\bigr)
-
\bigl((\partial_{\xi_\alpha}P^{(k)})P^{(j)}\bigr)\bigl(P^{(j)}\partial_{x^\alpha}P^{(l)}\bigr)
\bigr].
\end{multline*}
Using \eqref{dPP}, we can rewrite the above formula as
\begin{multline*}
{\operatorname{tr}}
\{P^{(k)},P^{(j)},P^{(l)}\}
={\operatorname{tr}}
\bigl[
\bigl(
\delta^{kj}\partial_{x^\alpha}P^{(j)}-P^{(k)}\partial_{x^\alpha}P^{(j)}
\bigr)
\bigl(
\delta^{jl}\partial_{\xi_\alpha}P^{(j)}-(\partial_{\xi_\alpha}P^{(j)})P^{(l)}
\bigr)
\\
-\bigl(
\delta^{kj}\partial_{\xi_\alpha}P^{(j)}-P^{(k)}\partial_{\xi_\alpha}P^{(j)}
\bigr)
\bigl(
\delta^{jl}\partial_{x^\alpha}P^{(j)}-(\partial_{x^\alpha}P^{(j)})P^{(l)}
\bigr)
\bigr].
\end{multline*}
Expanding the parentheses in the above formula and rearranging terms, we get
\begin{multline}
\label{trPkjl1}
{\operatorname{tr}}
\{P^{(k)},P^{(j)},P^{(l)}\}
=\delta^{kj}{\operatorname{tr}}\{P^{(j)},P^{(l)},P^{(j)}\}
+\delta^{jl}{\operatorname{tr}}\{P^{(j)},P^{(k)},P^{(j)}\}
\\
-\delta^{kl}{\operatorname{tr}}\{P^{(j)},P^{(k)},P^{(j)}\}.
\end{multline}
In the special case $l=k$ the above formula becomes
\begin{equation}
\label{trPkjk}
{\operatorname{tr}}\{P^{(k)},P^{(j)},P^{(k)}\}
=2\delta^{kj}{\operatorname{tr}}\{P^{(j)},P^{(j)},P^{(j)}\}
-{\operatorname{tr}}\{P^{(j)},P^{(k)},P^{(j)}\}.
\end{equation}
Each of the three terms in the RHS of
\eqref{trPkjl1}
can now be rewritten using the identity
\eqref{trPkjk} with appropriate choice of indices,
which gives us
\eqref{20 December 2016 equation 1}.

Let us now substitute
\eqref{20 December 2016 equation 1}
into the triple sum in the RHS of
\eqref{The matrix trace of the resolvent equation 3}:
\begin{multline}
\label{very long formula}
\sum_{j,k,l}
\frac{h^{(j)}-z}{(h^{(k)}-z)(h^{(l)}-z)}
{\operatorname{tr}}\{P^{(k)},P^{(j)},P^{(l)}\}
\\
=2\sum_{j}
\frac{1}{h^{(j)}-z}
{\operatorname{tr}}\{P^{(j)},P^{(j)},P^{(j)}\}
\\
-
\sum_{j,l}
\frac{1}{h^{(l)}-z}
{\operatorname{tr}}\{P^{(l)},P^{(j)},P^{(l)}\}
-
\sum_{j,k}
\frac{1}{h^{(k)}-z}
{\operatorname{tr}}\{P^{(k)},P^{(j)},P^{(k)}\}
\\
+
\sum_{j,k}
\frac{h^{(j)}-z}{(h^{(k)}-z)^2}
{\operatorname{tr}}\{P^{(k)},P^{(j)},P^{(k)}\}
\\
=2\sum_{j}
\frac{1}{h^{(j)}-z}
{\operatorname{tr}}\{P^{(j)},P^{(j)},P^{(j)}\}
-
2
\sum_{j,k}
\frac{1}{h^{(k)}-z}
{\operatorname{tr}}\{P^{(k)},P^{(j)},P^{(k)}\}
\\
+
\sum_{j,k}
\frac{h^{(j)}-z}{(h^{(k)}-z)^2}
{\operatorname{tr}}\{P^{(k)},P^{(j)},P^{(k)}\}
\\
=2\sum_{j}
\frac{1}{h^{(j)}-z}
{\operatorname{tr}}\{P^{(j)},P^{(j)},P^{(j)}\}
-
2
\sum_{k}
\frac{1}{h^{(k)}-z}
{\operatorname{tr}}\{P^{(k)},P^{(k)}\}
\\
+
\sum_{j,k}
\frac{1}{(h^{(k)}-z)^2}
{\operatorname{tr}}\{P^{(k)},h^{(j)}P^{(j)},P^{(k)}\}
-z
\sum_{k}
\frac{1}{(h^{(k)}-z)^2}
{\operatorname{tr}}\{P^{(k)},P^{(k)}\}
\\
=2\sum_{j}
\frac{1}{h^{(j)}-z}
{\operatorname{tr}}\{P^{(j)},P^{(j)},P^{(j)}\}
+
\sum_{k}
\frac{1}{(h^{(k)}-z)^2}
{\operatorname{tr}}\{P^{(k)},A_1,P^{(k)}\}
\\
=2\sum_{j}
\frac{1}{h^{(j)}-z}
{\operatorname{tr}}\{P^{(j)},P^{(j)},P^{(j)}\}
+
\sum_{j}
\frac{1}{(h^{(j)}-z)^2}
{\operatorname{tr}}\{P^{(j)},A_1,P^{(j)}\}
\\
=2\sum_{j}
\frac{1}{h^{(j)}-z}
{\operatorname{tr}}\{P^{(j)},P^{(j)},P^{(j)}\}
+
\sum_{j}
\frac{1}{(h^{(j)}-z)^2}
{\operatorname{tr}}\{P^{(j)},A_1-h^{(j)}I,P^{(j)}\},
\end{multline}
where we used the identities $\sum_jP^{(j)}=I$, $\{P^{(k)},P^{(k)}\}=0$
and \eqref{12 December 2016 equation 3}.
Substituting
\eqref{very long formula}
into
\eqref{The matrix trace of the resolvent equation 3}
we arrive at
\eqref{20 December 2016 equation 2}.

\end{appendices}


\begin{thebibliography}{19}

\bibitem{jst_review}
Z.~Avetisyan, Y.-L.~Fang and D.~Vassiliev,
Spectral asymptotics for first order systems.
{\it Journal of Spectral Theory} \textbf{6} (2016), 695--715.

\bibitem{jst_part_a}
O.~Chervova, R.~J.~Downes and D.~Vassiliev,
The spectral function of a first order elliptic system.
{\it Journal of Spectral Theory} \textbf{3} (2013), 317--360.

\bibitem{DuGu75}
J.~J.~Duistermaat and  V.~W.~Guillemin,
The spectrum of positive elliptic operators and periodic bicharacteristics.
{\it Invent.~Math.}~\textbf{29} (1975), 39--79.

\bibitem{DuiHor}
J.~J.~Duistermaat and L.~H\"ormander,
Fourier integral operators II.
{\it Acta~Math.}~\textbf{128} (1972), 183--269.


\bibitem{GrSj94}
A.~Grigis and J.~Sj\"ostrand,
{\it Microlocal analysis for differential operators, an introduction.}
London Math.~Soc.~Lect.~Notes ser.~196,
Cambridge University Press, 1994.

\bibitem{GrSe95}
G.~Grubb and R.~Seeley,
Weakly parametric pseudodifferential operators and Atiyah-Patodi-Singer boundary problems.
{\it Invent.~Math.}~\textbf{121} (1995), 481--529.

\bibitem{HeSj89}
B.~Helffer and J.~Sj\"ostrand,
 Semiclassical analysis for Harper's equation. III:  Cantor structure of the spectrum.
{\it Bull.~de la SMF} \textbf{117}(4)(1989), m{\'e}moire no 39. 

\bibitem{Ho85}
L.~H\"ormander, {\it The analysis of linear partial differential operators III. Pseudo-differential operators.} Grundlehren der Mathematischen Wissenschaften \textbf{274}. Springer-Verlag, Berlin, 1985. 

\bibitem{IvriiDoklady1980}
V.~Ivrii,
On the second term of the spectral asymptotics for the
Laplace--Beltrami operator on manifolds with boundary and
for elliptic operators acting in fiberings.
{\it Soviet Mathematics Doklady} \textbf{21} (1980), 300--302.

\bibitem{IvriiFuncAn1982}
V.~Ivrii,
Accurate spectral asymptotics for elliptic operators that act in vector bundles.
{\it Functional Analysis and Its Applications} \textbf{16} (1982), 101--108.

\bibitem{ivrii_springer_lecture_notes}
V.~Ivrii,
{\it Precise spectral asymptotics for elliptic operators
acting in fiberings over manifolds with boundary.}
Lecture Notes in Mathematics \textbf{1100}, Springer-Verlag, Berlin, 1984.

\bibitem{ivrii_book}
V.~Ivrii,
{\it Microlocal analysis and precise spectral asymptotics.}
Springer-Verlag, Berlin, 1998.

\bibitem{kamotski}
I.~Kamotski and M.~Ruzhansky,
Regularity properties, representation of solutions, and spectral asymptotics
of systems with multiplicities.
{\it Comm.~Partial Differential Equations} \textbf{32} (2007), 1--35.

\bibitem{mckeag}
P.~McKeag and Yu.~Safarov,
Pseudodifferential operators on manifolds: a coordinate-free approach.
In Partial Differential Equations and Spectral Theory,
part of the Operator Theory: Advances and Applications book series, volume \textbf{211},
Springer, Basel, 2011, 321--341.

\bibitem{grisha}
G.~V.~Rozenblyum,
Spectral asymptotic behavior of elliptic systems.
{\it Journal of Mathematical Sciences} \textbf{21} (1983), 837--850.


\bibitem{SafarovDSc}
Yu.~Safarov,
Non-classical two-term spectral asymptotics for self-adjoint elliptic operators.
DSc~thesis,
Leningrad Branch of the Steklov Mathematical Institute of the
USSR Academy of Sciences, 1989.
In Russian.

\bibitem{mybook}
Yu.~Safarov and D.~Vassiliev,
{\it The asymptotic distribution of eigenvalues of partial differential operators.}
Amer.~Math.~Soc., Providence (RI), 1997.

\bibitem{shubin}
M.~A.~Shubin,
{\it Pseudodifferential operators and spectral theory.}
Springer-Verlag, Berlin, second edition, 2001.

\bibitem{SjZw02}
J.~Sj\"ostrand and M.~Zworski, 
Quantum monodromy and semi-classical trace formulae.
{\it Journal de Math{\'e}matiques Pures et Appliqu{\'e}es} \textbf{81}(2002), 1--33. 


\end{thebibliography}
\end{document}